\documentclass{article}

\usepackage{PRIMEarxiv}
\usepackage{cite}
\usepackage[utf8]{inputenc} %
\usepackage[T1]{fontenc}    %
\usepackage{hyperref}       %
\usepackage{url}            %
\usepackage{booktabs}       %
\usepackage{amsfonts}       %
\usepackage{nicefrac}       %
\usepackage{microtype}      %
\usepackage{lipsum}
\usepackage{fancyhdr}       %
\usepackage{graphicx}       %
\graphicspath{{Figures/}}     %
\usepackage{amsmath}

\title{Dynamical System Parameter Path Optimization using Persistent Homology}

\author{
  Max M. Chumley \\
  Mechanical Engineering and \\Computational Mathematics, Science and Engineering \\
  Michigan State University \\
  East Lansing MI\\
  \texttt{chumleym@msu.edu} \\
   \And
  Firas A. Khasawneh \\
  Computational Mathematics, Science and Engineering \\ 
  Michigan State University \\
  East Lansing MI\\
  \texttt{khasawn3@msu.edu} \\
}

\begin{document}
\maketitle

\begin{abstract}
    Nonlinear dynamical systems are complex and typically only simple systems can be analytically studied. In applications, these systems are usually defined with a set of tunable parameters and as the parameters are varied the system response undergoes significant topological changes or bifurcations. In a high dimensional parameter space, it is difficult to determine which direction to vary the system parameters to achieve a desired system response or state. In this paper, we introduce a new approach for optimally navigating a dynamical system parameter space that is rooted in topological data analysis. Specifically we use the differentiability of persistence diagrams to define a topological language for intuitively promoting or deterring different topological features in the state space response of a dynamical system and use gradient descent to optimally move from one point in the parameter space to another. The end result is a path in this space that guides the system to a set of parameters that yield the desired topological features defined by the loss function. We show a number of examples by applying the methods to different dynamical systems and scenarios to demonstrate how to promote different features and how to choose the hyperparameters to achieve different outcomes.  
\end{abstract}

\keywords{Dynamical Systems \and Bifurcations \and Persistent Homology \and Topological Data Analysis \and Parameter Space Navigation}

\section{Introduction}

Dynamical systems are often defined in terms of tunable parameters and as the parameters are varied, the state or behavior of the system changes. These changes can be desired such as reducing oscillations in a mechanical system for safety or they can be detrimental such as a system changing from predictable periodic oscillations to chaotic behavior which is much more difficult to accurately predict. These changes in dynamic state are called bifurcations. It is assumed throughout that we have access to sampled realizations $X = [x_1,\cdots,x_N]$ where $x_i \in \mathbb{R}^n$ of a nonlinear dynamical system, and that $\vec{\mu}\in\mathbb{R}^D$ is the vector of system parameters. While for this work we generate $X$ from a known model $\dot{x}=f(x,t, \vec{\mu})$, in theory $X$ could be generated from experimental data, but many experiments need to be conducted to adequately sample the parameter space. Bifurcations can occur in both continuous and discrete time dynamical system and they are characterized by qualitative changes in the response as one or more parameters (called the bifurcation parameters) are varied. Bifurcations have been studied extensively in the literature for $D=1$ \cite{Kantz2004,Hirsch2003,Baker1996,Sayama2015}, however, when $D>1$, it is much more difficult to study how the system response varies in the parameter space. This often leads to researchers fixing parameters to reduce the dimensionality of $\vec{\mu}$ which only yields information on a small projection of the parameter space. Bifurcations typically indicate that the system is transitioning from one state to a topologically differing state in the phase space. One visual tool for finding bifurcations is the bifurcation diagram ($D=1$), which shows local extrema of a given system over a varying control parameter while keeping other parameters fixed. As more bifurcation parameters are added, locating bifurcations becomes more difficult as infinitely many paths exist between any two points in the parameter space. In this work, we aim to locate optimal paths in this space that connect an initial state to a desirable state while avoiding regions of the space that lead to unsafe or undesired dynamics. We assume that system parameters cannot be changed discontinuously or in other words we cannot teleport from one point in the parameter space to another and we only have information that is local to the current $\vec{\mu}$.

When the governing equations for a deterministic dynamical system are available, then there are tools that facilitate tracking the bifurcations as a parameter varies; although, exhaustively tracking all the bifurcations is not a trivial task. 
One such tool is numerical continuation \cite{Kuznetsov1998,Seydel2009,Dankowicz2013}, which is a path following approach that tracks the solution branches as system parameters are varied. 
However, if the governing equations are not available or are too complicated, then sometimes it is possible to track the solutions and the bifurcations of the underlying dynamical system using Control-Based Continuation (CBC) \cite{Sieber2007,Sieber2010,Barton2011,Bureau2013,Barton2013,Barton2017}. 
CBC was successfully used in many scientific domains including biochemistry \cite{Godwin2017}, physics \cite{KRAUSKOPF2006}, mechanics \cite{Peeters2009}, and fluid dynamics \cite{Huntley2017}. 
In this setting, numerical continuation is applied to a feedback-controlled physical experiment such that the control becomes non-invasive \cite{Barton2017}. Treating the physical system as a numerical model, control-based continuation allows systematic investigations of the bifurcations in the system by treating the control target as a proxy for the state. 
However, due to the stabilizing feedback, stability information and consequently the ability to detect and classify certain bifurcations is lost. 
Nevertheless, existing tools for tracking bifurcation or exploring dynamic changes in state space remain limited to small spaces with most of the time one and at most two bifurcation control parameters. High-dimensional parameter spaces of dynamical systems have been explored using features of numerical solutions in \cite{kuehn2008}, but this method relies on choosing a feature of the system response which is unintuitive and it also uses random exploration and interpolation of the parameter space with interpolation to extract information from the system. 

Therefore, there is a need for an intuitive, data-driven approach to navigating high dimensional dynamical system parameter spaces to guide the system to an acceptable response. We set out to develop a framework with topological data analysis and persistence optimization at its core to meet this need and provide a more intuitive understanding of the map between parameter space dynamics and topological persistence. We accomplished this goal in three phases. First, we 
provide the necessary background for the tools we used from topological data analysis and 
persistence optimization followed by a dictionary of cost function terms to promote different persistence diagrams that map to dynamical system responses in Section~\ref{sec:cost_lib}. For the second phase, we show preliminary derivative-free optimization results moving from chaotic behavior to periodic in the Lorenz system in Section~\ref{sec:preliminary_path}. For phase 3, we perform the optimization using gradient descent with the cost function library from phase 1 in Section~\ref{sec:parameter_nav} and show extensive numerical results in \ref{sec:param_opt_experiments}.

\section{Methods}

Topological Data Analysis (TDA) is a collection of tools that can generally be used to quantify shape or structural information from data in different formats. The flagship tool from TDA is called Persistent Homology (PH) and specifically in this work we focus on PH computed on point clouds in $\mathbb{R}^n$. In Section~\ref{sec:ph} we provide a basic description and example of PH and in Section~\ref{sec:ph_opt} we show how functions of persistence can be differentiated to enable gradient descent for optimizing topological features of data. More specifics on PH can be found in \cite{Hatcher2002,Kaczynski2004,Ghrist2008,Carlsson2009,Edelsbrunner2010,Mischaikow2013,oudot2017persistence,Munch2017} for the interested reader. In Section~\ref{sec:parameter_nav}, we show how the differentiability pipeline is augmented with an ODE solver and how the map between the parameter space and the state space point cloud is differentiated to enable gradient descent through the full map. Lastly, we define the topological cost function library in Section~\ref{sec:cost_lib} by mapping different features in the persistence diagram to dynamical system behaviors.

\subsection{Persistent Homology}\label{sec:ph}

\begin{figure}[htbp]
    \centering
    \includegraphics[width=\textwidth]{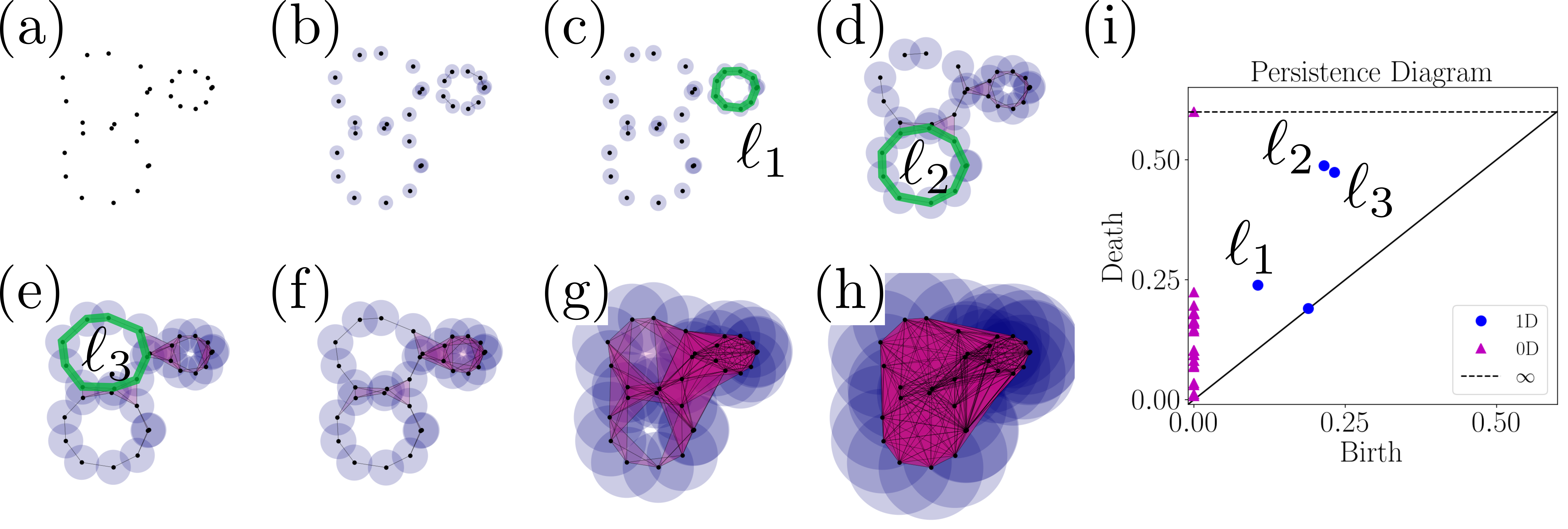}
    \caption{Point cloud persistent homology example. (a-h) show the Vietoris-Rips filtration for increasing values of the connectivity parameter and (i) shows the full persistence diagram with the 3 prominent loops labeled in the persistence diagram and the filtration.}
    \label{fig:ph_example}
\end{figure}

With persistent homology, we start by inducing a simplicial complex or generalized graph on the data based on a set of rules that are defined by a chosen filter function with a varying connectivity parameter. As the connectivity parameter is varied, simplices (edges, faces, etc.) are added to the simplicial complex such that each successive complex includes the previous. This process is called a filtration. For example, in Fig.~\ref{fig:ph_example}(a) we see the starting point cloud with no edges or faces and in Fig.~\ref{fig:ph_example}(b) after the connectivity parameter increases the original complex in Fig.~\ref{fig:ph_example}(a) is contained in the complex in Fig.~\ref{fig:ph_example}(b). In this example, the connectivity parameter, $r$, is represented as the radius of balls that are centered at each point in the point cloud and when the balls intersect, an edge is added between those vertices and for any 3 balls intersecting, a face or triangle is added. This specific simplicial complex is called the Vietoris-Rips (VR) complex. For a given value of $r$, the induced simplicial complex has inherent topological features or homology such as connected components (0D homology) or holes (1D homology). For example, in Fig.~\ref{fig:ph_example}(c) we see the loop $\ell_1$ has formed or is \textit{born} in the simplicial complex and as $r$ increases, $\ell_1$ eventually fills in or \textit{dies} in Fig.~\ref{fig:ph_example}(f). We track the birth and death of homological features as $r$ changes in a persistence diagram Fig.~\ref{fig:ph_example}(i) where the birth and death of topological features are plotted as $x$ and $y$ coordinates. Persistence diagrams have been proven to be stable under small perturbations of the point cloud \cite{CohenSteiner2006} and provide a compact representation of complex topological features of the data. 

\subsection{Persistence Optimization}\label{sec:ph_opt}

\begin{figure}[p]
    \centering
    \includegraphics[width=\textwidth]{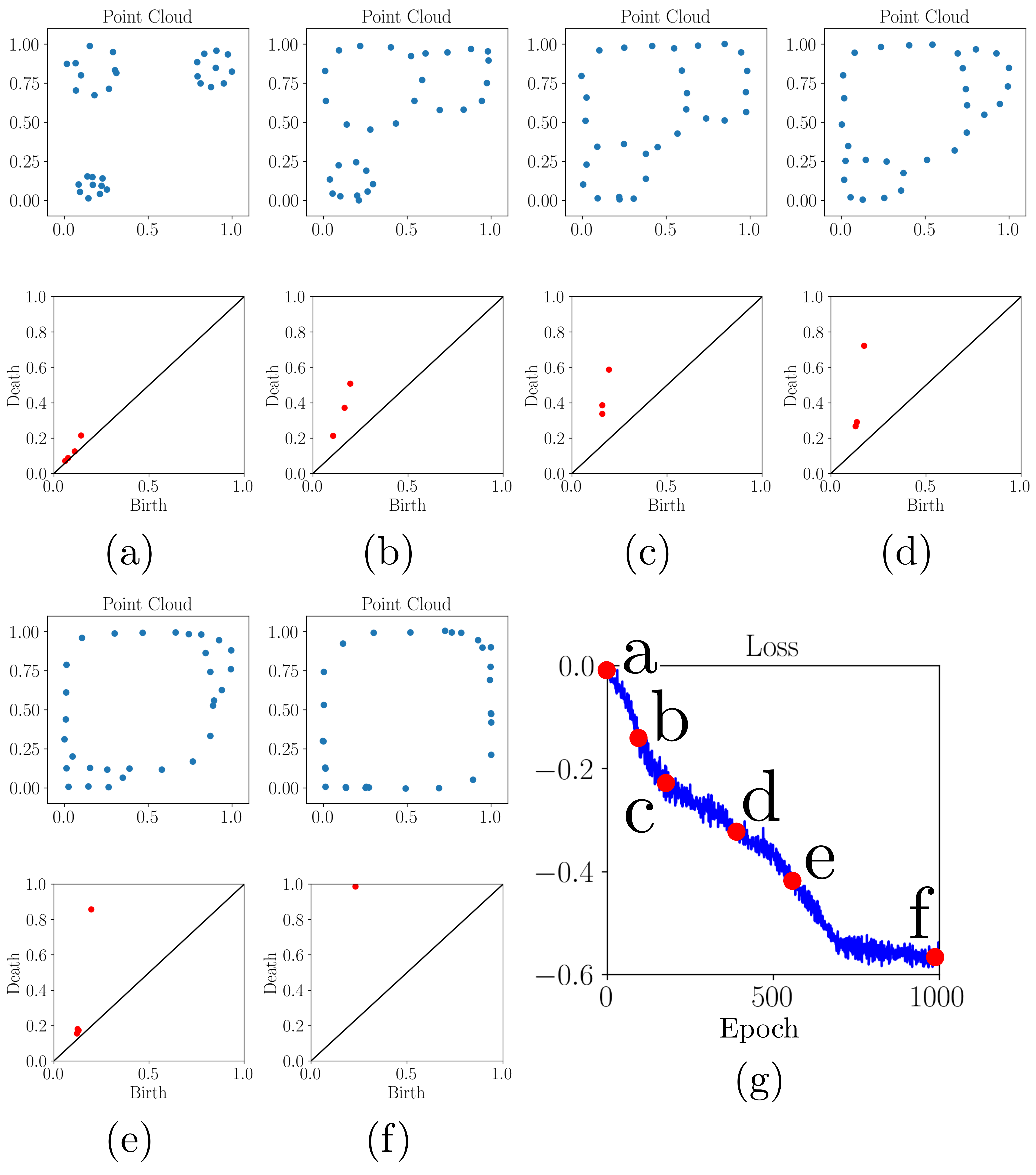}
    \caption{Persistence optimization example. The initial point cloud is shown in (a)  and the updated point cloud and persistence diagrams are shown at different points in the optimization process in (b-f). The loss is plotted in (g) with the different points labeled at the respective epoch numbers.}
    \label{fig:pd_opt_example}
\end{figure}

In more recent literature, the differentiability of persistence diagrams has been defined allowing for integrating persistent homology into an optimization pipeline \cite{Carriere2020,Leygonie2021}. This pipeline is represented by $\mathcal{M} \xrightarrow[\hspace{1.25cm}]{B} \text{Pers} \xrightarrow[\hspace{1.25cm}]{V} \mathbb{R}$ where $\mathcal{M}$ is the point cloud and $B$ maps the point cloud to its persistence diagram. $V$ is the function of persistence that maps to a real number. Using the chain rule, the map $V\circ B$ is differentiated to perform gradient descent and perturb $\mathcal{M}$ to minimize $V\circ B$. Specifically, the derivative of a persistence map $B$ with respect to a perturbed persistence map $\tilde{B}$ and the point cloud $X$ using the Vietoris-Rips complex is computed using,
\begin{equation}\label{eq:pd_der}
    d_{P,\tilde{B}}B(\hat{u})=\left[\left(P^T_{v(\sigma),w(\sigma)}\hat{u},~P^T_{v(\sigma '),w(\sigma ')}\hat{u}  \right)_{i=1}^m\right],
\end{equation}
where $\hat{u}$ is a chosen perturbation of the point cloud, $P$ tracks the changes in attaching edges for each persistence pair in the persistence diagram or the edge $\sigma$ that result in the birth of the $i$-th persistence pair and the face $\sigma '$ that results in the death of that same persistence pair. $v$ and $w$ are the vertices of attaching edges for $\sigma$ or $\sigma '$\cite{Leygonie2021}. $P$ is computed using $P_{i,j}=\frac{p_i-p_j}{||p_i-p_j||_2}$ where $p_i$ corresponds to the $i$-th attaching edge vertex. As a result, the derivative of a persistence diagram gives a set of vectors with one for each persistence pair which quantifies how the persistence diagram changes for a particular perturbation $\hat{u}$ \cite{Leygonie2021}. A helpful way of thinking about this process is to imagine tracking the edges that result in the birth of a persistence feature and the simplex that results in the death of that feature and quantify how much those birth and death values change if the input point cloud is infinitesimally changed in some direction. In \cite{Carriere2020}, the authors detail how this optimization process works and what conditions need to be met. Specifically the point cloud needs to be in general position for the derivative to be defined so no points can be in the same position and no pairs of points can be equidistant. The \texttt{gudhi} Python library supports \texttt{PyTorch} auto-differentiation for performing this optimization \cite{Carriere2020}. An example of this process is shown in Fig.~\ref{fig:pd_opt_example} where the loss function is defined to maximize the sum of the 1D persistence lifetimes and regularization is applied to keep the points in the range 0-1. Figure~\ref{fig:pd_opt_example}(a) shows the original point cloud and persistence diagram defined as three separate loops. As the optimization is applied using gradient descent, the loops grow in Fig.~\ref{fig:pd_opt_example}(b-e) and eventually merge into a single loop in Fig.~\ref{fig:pd_opt_example}(f) to minimize the loss function. Finally, the loss plot with respect to the epoch is shown in Fig.~\ref{fig:pd_opt_example}(g). We see that the loss reaches a minimum of around -0.6 after 1000 epochs. This framework enables optimizing topological features of the data that would otherwise be difficult to describe and due to the tight link between topology and dynamical systems it is an ideal candidate for defining an intuitive topological language for moving through the parameter space of a dynamical systems.

\subsection{Gradient Descent Parameter Space Navigation using TDA}\label{sec:parameter_nav}

\begin{figure}[htbp]
  \centering
  \includegraphics[width=\textwidth]{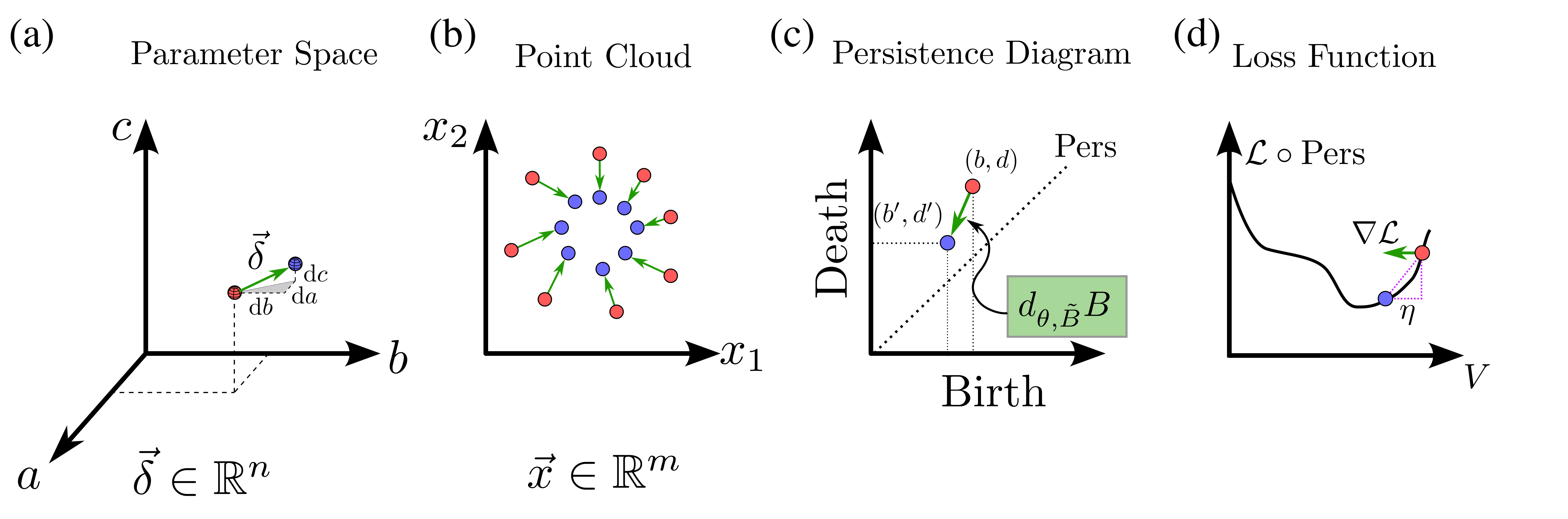}
  \caption{Diagram demonstrating the map from the parameter space to the loss function as solving the inverse problem from taking a step against the gradient of the loss function to reach a new set of parameters in the parameter space propagated through the persistence diagram and the state space point cloud.}
  \label{fig:par_diff}
\end{figure}

The main goal of this work is to augment this pipeline with the parameter space of a dynamical system. Specifically, the map becomes $\mathcal{D} \xrightarrow[\hspace{1.25cm}]{B'} \mathcal{M} \xrightarrow[\hspace{1.25cm}]{B} \text{Pers} \xrightarrow[\hspace{1.25cm}]{V} \mathbb{R}$, where $\mathcal{D}$ is the parameter space of a dynamical system and $B'$ is the map from the parameter space to the state space point cloud. The rest of the pipeline remains the same because $\mathcal{M}$ is still a point cloud, but its movement is governed by the parameter space and dynamics of the system. To have differentiability in the original pipeline from \cite{Leygonie2021}, the authors placed restrictions on the positions of point to ensure the derivative exists. Specifically the point cloud must be in general position. Similar restrictions need to be considered for $B'$ before this augmented pipeline can be used. $B'$ is essentially the numerical integrator for the system that generates the state space of the system. Differentiation of ODE solvers has been enabled using the \textit{adjoint sensitivity method} \cite{chen2019} where the gradient of the loss function is computed for each state of the system by solving another ODE whose solution gives the gradient of each state of the system trajectory. This method has been implemented in many different solvers in the torchdiffeq python library with pytorch compatibility \cite{chen2019}. Typically this method is used for neural ODEs which are a continuous analog of traditional neural networks, but for this work we only need the ability to compute the gradient of $B'$. The method works by assuming we have a dynamical system $\dot{x}=f(x(t), t, \mu)$ with some loss function that depends on the states of the system $x(t)$ given by,
\begin{equation}
  L(x(t))=L \left(x(t_0) + \int_{t_0}^{t_f} f(x(t),t,\mu)dt\right).
\end{equation}
The goal is to obtain the gradient of $L$ with respect to the system parameters $\mu$. In \cite{chen2019}, the authors define an adjoint state $a(t)=\partial L/\partial x(t)$ and show a different ODE that governs its dynamics. From \cite{chen2019}, the gradient of the loss function with respect to $\mu$ is then given by the integral,
\begin{equation}
  \frac{\partial L}{\partial \mu} = -\int_{t_0}^{t_f} a(t)^T \frac{\partial f(x(t),t,\mu)}{\partial \mu}dt,
\end{equation} 
where $a(t)$ is obtained by solving,
\begin{equation}
  \frac{d a(t)}{dt} = -a(t)^T\frac{\partial f(x(t),t,\mu)}{\partial x}.
\end{equation}
We see that ultimately, the gradient of the loss with respect to system parameters depends on how much the states $x(t)$ change with respect to changes in the system parameters so if a slight change in parameters yields drastically different system states, $\partial L/\partial \mu$ can be very large and cause complications in the optimization process. However, this method shows that it is possible to compute the gradient of the map $B'$ in our pipeline allowing for the full inverse problem to be solved as shown in Fig.~\ref{fig:par_diff} where a 3-dimensional parameter space is shown in Fig.~\ref{fig:par_diff}~(a) with starting point in red and to move to the next path point, a step is taken toward a minimizer of the loss function in Fig.~\ref{fig:par_diff}~(d) and the step is propagated back using gradient descent through the persistence diagram (Fig.~\ref{fig:par_diff}~(c)) and state space (Fig.~\ref{fig:par_diff}~(b)) to obtain a direction in the parameter space for the next point on the path. This enables a gradient descent approach to moving through the parameter space using the full persistence diagrams to move to a set of parameters where the response meets the criteria defined from Section~\ref{sec:cost_lib}. As a result, the gradient of the full map $V\circ B\circ B'$ is computed to move through the parameter space leveraging the differentiability of persistence diagrams \cite{Gameiro2016, Carriere2020, Leygonie2021} and the adjoint sensitivity method from \cite{chen2019}.

\subsection{Cost Function Library}\label{sec:cost_lib}

Currently, there is only a basic understanding of the general shape of a persistence diagram for a given dynamic state. For example a periodic response often contains a single 1D persistence pair with a long lifetime. We aim to create a dictionary of persistence diagrams with different traits that will allow the user to impose constraints on the problem. By combining these criterion, a \textit{desired} persistence diagram will be obtained effectively designing an objective function for the optimization problem using topological characteristics of acceptable system behavior. So in this setting, we assume that there is a target persistence diagram that corresponds to desirable criteria for the system response. In this work, we focus on the ability to promote or deter the following behaviors: periodic solutions, fixed point stability, chaotic behavior and allowing for specifying regions of the parameter space that are off limits or unsafe for the system. 

\subsubsection{Periodic Solutions}

To promote periodic solutions, it is intuitive to see that the persistence diagram should contain at least one long lifetime pair far from the diagonal. Therefore, the maximum 1D persistence lifetime feature can be used to control the size of the largest loop. Maximizing the maximum persistence encourages larger loops in the state space. This feature is computed using,
\begin{equation}
  \text{maxPers}_1 = \max_i{(\ell_i^1)},
\end{equation}
where $\ell_i^1$ is the lifetime of the $i$th 1D loop in the persistence diagram. Another feature that is typically used to quantify the prominence of persistence pairs is the total persistence where instead of taking the maximum lifetime, the sum of all lifetimes is used. However, to promote periodic solutions, we argue that the maximum lifetime is more important than the total lifetime because if all lifetimes are maximized simultaneously this could also promote chaotic behavior. The maximum persistence could also promote chaos so it needs to be combined with other cost function terms such as persistent entropy if there is chaotic behavior in the system. 

\subsubsection{Fixed Point Stability}

The second criterion is to promote fixed point stability. In this case, the persistence diagram should contain loops that are close to the diagonal. This is because fixed points are typically represented by points in the state space that are close in proximity. To quantify this behavior, the maximum persistence can also be used and minimizing this term will result in loops that are low lifetime. Another feature that could be used is the total persistence feature of the 1D persistence diagram to promote all loop lifetimes approaching zero. This feature is computed as,
\begin{equation}
  \text{totPers}_1 = \sum_i{\ell_i^1}.
\end{equation}
The average persistence can also be used which is the total persistence divided by the number of persistence pairs to normalize the feature. In some cases these features may be biased by persistence pairs near the diagonal so in this case maximum persistence should be used or the $n$ longest lifetime persistence pairs to filter out low lifetime features.

\subsubsection{Chaos}

By definition, chaos in a dynamical system is a sensitive dependence of initial conditions or parameters. In other words, changing a parameter or starting point by an infinitesimal amount will yield drastically different system trajectories over time, but the topologies of the trajectories should be similar. A 1D point cloud persistence diagram for a chaotic trajectory typically consists of many loops that are close and moderately far from the diagonal. To control chaos using persistence diagrams, we suggest using persistent entropy which is the Shannon entropy for a probability distribution of persistence lifetimes \cite{Atienza2020}. Specifically, persistent entropy is computed as,
\begin{equation}
  E= -\sum_i {p_i\log_2{(p_i)}},
\end{equation}
where $\ell^1$ is the set of 1D persistence lifetimes, $p_i=\ell_i^1/L$ and $L=\sum_i {\ell_i^1}$. Persistent entropy gives a measure of order for the distribution of persistence lifetimes so if the lifetimes are more unevenly distributed, persistent entropy is larger and for lifetimes that are more concentrated, $E$ is smaller. Persistent entropy was shown to be stable under small perturbations in \cite{Atienza2020} and is a Lipschitz function so it can be used with persistence optimization. Note that persistent entropy is biased by the number of persistence pairs in the persistence diagram so to reduce this effect it is often normalized by $\log_2{(N)}$ where $N$ is the number of pairs in the persistence diagram. This ensures that only the level of disorder is being quantified. While in theory this function can be used for quantifying disorder in persistence diagrams and allow for moving to parameters with a more ordered lifetime distribution, we will see in Section~\ref{sec:parameter_nav} that there are complications in computing the gradients of a system trajectory in a chaotic region of the parameter space that need to be mitigated using different optimization techniques. Examples using this term to define cost functions are shown in Section~\ref{sec:param_opt_experiments}.

\subsubsection{Forbidden Regions}\label{sec:forbidden}

The last cost function term we consider deals with allowing the user to specify regions of the parameter space that are forbidden. In other words, if the parameters enter these regions, a penalty is applied to the cost function. We assume that there is a function $f(\vec{\mu}): \mathbb{R}^m \to \mathbb{R}$ that is positive for values of $\vec{\mu}$ that are inside of the forbidden region bounded by $f$ and negative for $\vec{\mu}$ outside of this region. Here, $\vec{\mu}$ is a vector of system parameters and $m$ is the number of system parameters or the dimension of the parameter space. A penalty term $\mathcal{L}_p$ can be added to the overall cost function in the form of,
\begin{equation}\label{eq:forbidden}
  \mathcal{L}_p = \exp{(a~f(\vec{\mu}))}
\end{equation}
where $a\in\mathbb{R}^+$ is a parameter that is chosen based on tolerances for how close the parameters are allowed to be to the boundary. For example, if $a$ is close to zero, $\mu$ will be strongly forced away from $f$, but for large $a$ the parameters are allowed to get very close to the boundary and once it is crossed a large penalty is applied. $a$ needs to be balanced to ensure that gradients are not too discontinuous. For this work we use a value of $a=100$ unless otherwise specified. This logic can also be applied to the persistence pairs. For example, if we replace $\vec{\mu}$ with a persistence diagram and define $f$ to penalize persistence pairs in a forbidden region of the birth--death plane, the allowable areas for persistence pairs can also be specified. In all examples shown in Section~\ref{sec:param_opt_experiments}, penalty terms are used to regularize the paths by ensuring that the parameters stay in a specified region of the space. In this case, each $f$ is defined as a difference between each parameter in $\mu$ and its maximum and minimum allowable value. If the parameter leaves this range the penalty term increases. 

\begin{figure}[htbp]
  \centering
  \includegraphics[width=0.6\textwidth]{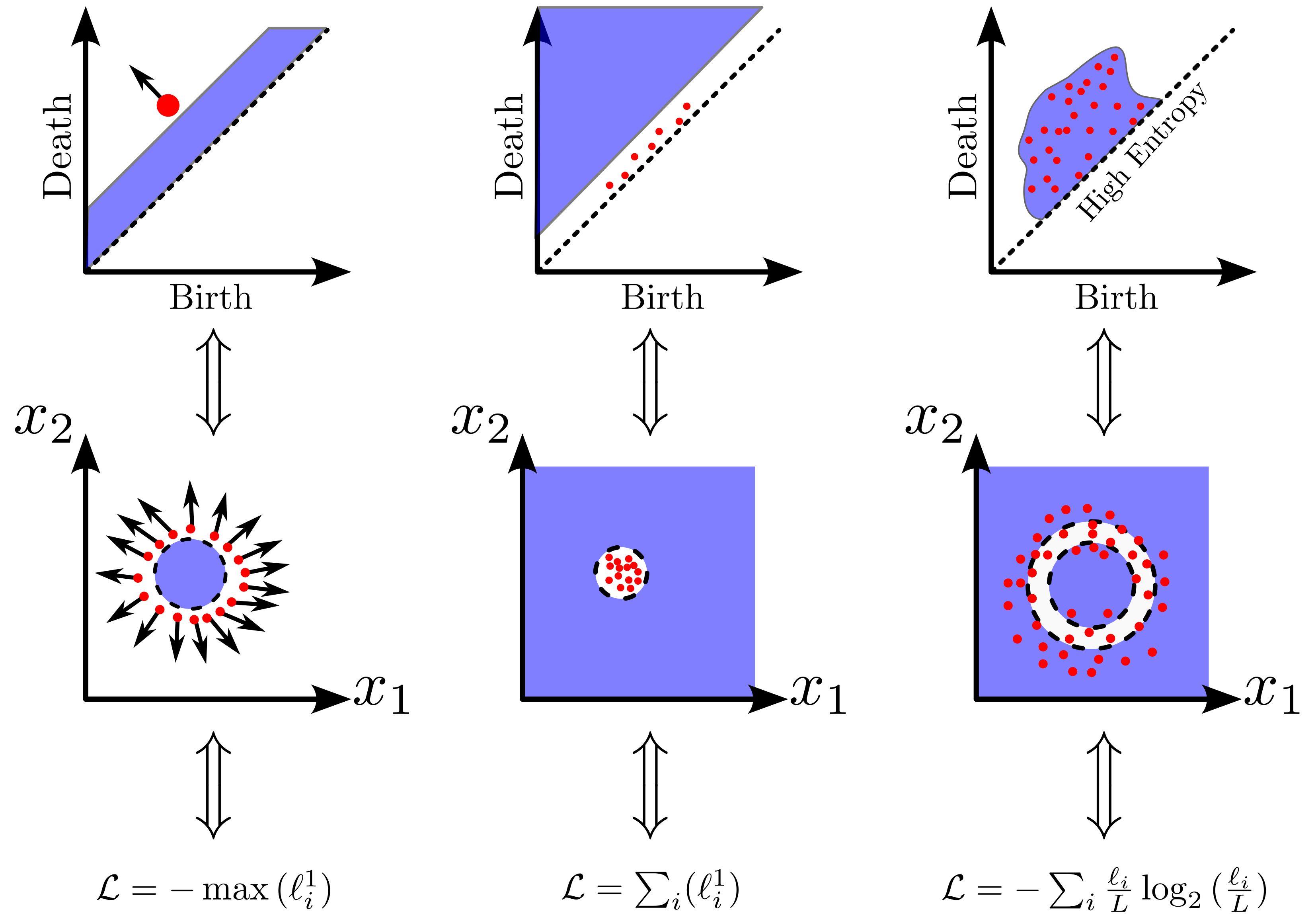}
  \caption{Example response criteria mapped into persistence diagrams. (Left) Maximizing maximum persistence to promote a large loop in the state space, (Middle) Limiting persistent loops to be close to the diagonal to encourage fixed point behavior, and (Right) Using high persistent entropy to classify chaotic regions in the persistence diagram.}
   \label{fig:pd_criteria}
\end{figure}

These criteria can be combined to build cost functions that promote desired responses. For example, the user may be searching for parameters that will result in a periodic solution with an amplitude less than a certain value. It is easy to see that this would map into the persistence diagram space as a 1D loop with a limited lifetime or death time minus birth time. In this case the maximum persistence can be used to promote a large loop as shown in Fig.~\ref{fig:pd_criteria}~(left) and can be combined with a penalty term from Eq.~\eqref{eq:forbidden} to deter lifetimes above a specified value which directly corresponds to the size of the loop or amplitude. A corresponding state space representation of this idea is shown in the middle row and potential cost function terms for achieving these behaviors are shown on the bottom row. 

Another desired behavior could be for the system to have fixed point stability. 
In this case the desired persistence diagram should have all of the loops close to the diagonal as shown in Fig.~\ref{fig:pd_criteria}~(middle) and the state space plot would have localized points to promote steady state stability. The third case shown in Fig.~\ref{fig:pd_criteria}~(right) is an avenue for classifying chaotic behavior using the persistence diagram by computing persistent entropy as in \cite{Myers2019}. In the state space this could correspond to a safe region being within a small annular region or a closed curve. Together these criteria are specified by the user form the desired characteristics of the target persistence diagram which are used for intuitive loss function engineering for computing the optimal path in the parameter space.   

\section{Preliminary Results}\label{sec:preliminary_path}

Once the desired persistence diagram is identified by defining a loss function to promote desired features, the objective is to move to a point in the parameter space that results in obtaining a response that has a persistence diagram closest to the desired diagram or in other words minimizes the loss function. Before gradient descent is performed with this cost function, we show preliminary results to further motivate the need for this approach and show that it is feasible using derivative free optimization methods.

\begin{figure}[thbp]
	\centering
	\includegraphics[width=0.8\textwidth]{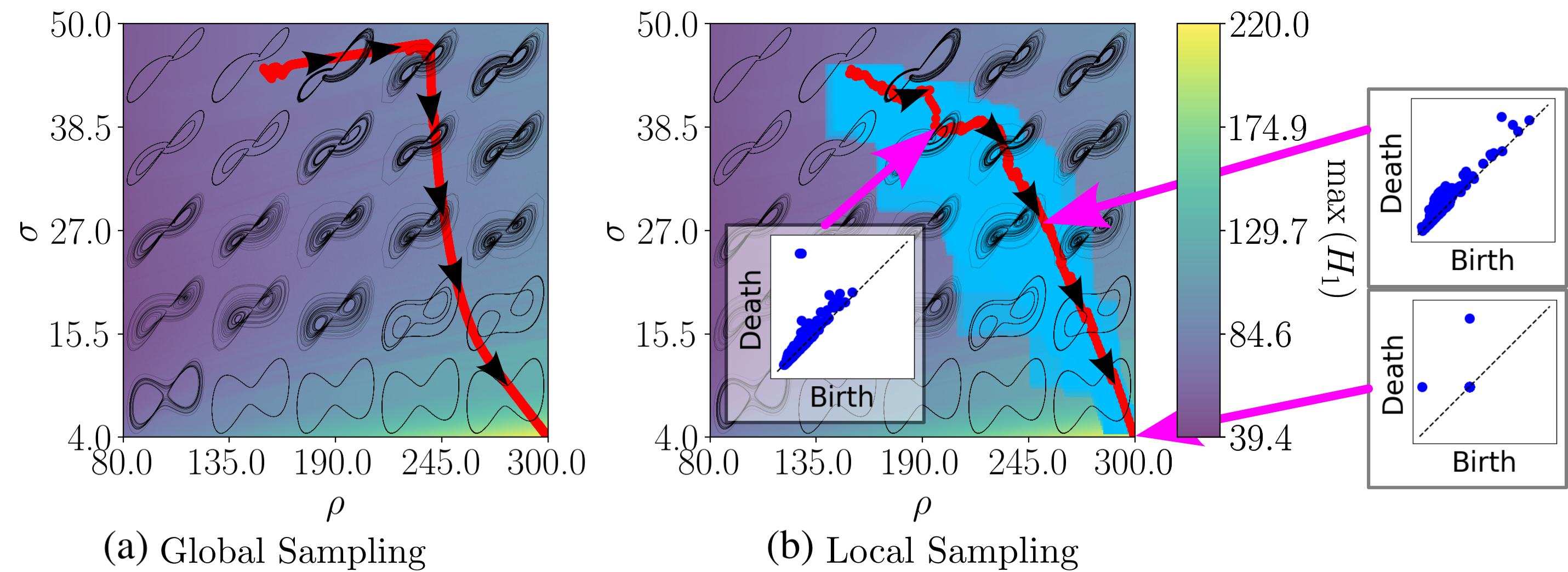}
	\caption{Lorenz system optimal parameter space paths using the global and local updating schemes. Corresponding persistence diagrams are shown at three points to demonstrate the topological differences between dynamic states. }
	\label{fig:lorenz_paths}
\end{figure}

Figure~\ref{fig:lorenz_paths} shows preliminary results for this work. Here, we considered the response optimal if it had the largest 1D persistence lifetime in the region of the parameter space being searched. This objective was constructed to promote periodic solutions over chaotic by maximizing the maximum persistence feature. Consider the Lorenz system given by $\dot{x}=\sigma(y-x),~\dot{y}=x(\rho-z)-y,~\dot{z}=xy-\beta z,$ where $\sigma$, $\rho$ and $\beta$ are system parameters and $x$, $y$, and $z$ are the system states. 
We restrict the parameter space to the plane $\beta=8/3$ for visualization and search for an optimal two dimensional path in the $(\rho,\sigma)$ space. We further limit the parameter space by setting $\rho\in[80,300]$ and $\sigma\in[4,50]$. The system was then simulated over a 500$\times$500 grid of parameters in the parameter space and the maximum 1D persistence lifetime was plotted as an image along with a subset of the system trajectories in two dimensions as shown in Fig.~\ref{fig:lorenz_paths}. 
For a given starting point in the parameter space, we aim to navigate this space optimally to reach the point with the largest 1D persistence lifetime. In this region of the parameter space, the optimal point was found to be $(\rho, \sigma)=(300,4)$ using the simplicial homology global optimization (SHGO) method from \cite{Endres2018}. 

\subsection{Navigation Schemes}
To generate the next point along the path towards the target state in Fig.~\ref{fig:lorenz_paths}, we used conventional global optimization algorithms to find the maximum 1D persistence in a local region near the starting point. A smaller sampling region highlighted in blue in Fig.~\ref{fig:lorenz_paths}b was chosen to obtain a smoother path to the optimal point. These two possible sampling schemes are described in the following sections.

\textbf{Global Sampling:} The first sampling method works by forming a rectangular region that grows from the starting point in the parameter space and solving for the global optimizer within that region. Let $(x_0, y_0)$ be the starting point in the 2D parameter space and the global problem domain $\Omega=\{\vec{\mu}=(x, y)\in[x_{min},x_{max}]\times[y_{min},y_{max}]\}$. The local search region is then generated as a fraction of the global region by the sequence, $\Omega_k=\{(x, y)\in\frac{1}{N}[(N-k)x_0+kx_{min},(N-k)x_0+kx_{max}]\times\frac{1}{N}[(N-k)y_0+ky_{min},(N-k)y_0+ky_{max}]: k=1...N\}$ where $N$ is the number of desired path steps. So as the step index $k$ increases, the feasible region grows to fill the entire global region when $k=N$. 
At each step $k$, we solve $\mu_{k}=\arg \max_{\mu\in\Omega_k} f(\vec{\mu})$ to find the direction vector relative to the current point. For this example, $f(\vec{\mu})=\max{(H_1)}$ is the maximum 1D persistence lifetime of the system simulation point cloud at a parameter input $\vec{\mu}$. Let $\hat{\mu}=\frac{\mu_{k}-\mu_{k-1}}{||\mu_{k}-\mu_{k-1}||}$ be the optimal unit direction from the local optimization problem assuming the optimal point is not identical to the current point. If this is the case, the path remains at the current point. If $\hat{\mu}$ is nonzero, the next point on the path is computed as $(x_k,y_k)=\alpha\hat{\mu}$ where $\alpha$ is the step size. 
Applying this updating scheme to the Lorenz system with a starting parameter vector of $(\rho, \sigma)=(153,45)$ with a constant step size of 0.1 and path length of 2500 steps, we obtained the path shown in Fig.~\ref{fig:lorenz_paths}~(a). It is clear that as more steps are taken in the path, the algorithm eventually moves in the direction of the optimal point and approaches it by the final step demonstrating optimal movement in the parameter space to move the system to a periodic response. 

\textbf{Local Sampling:} The global sampling path method required data from the full parameter space sampling region for solving the 2500 optimization problems. Simulation data is not always abundantly available so it is important to have an algorithm that also minimizes the search region size for the individual problems. To improve the path generation algorithm, we aim to use sampling regions that are centered around the current point essentially forming a rectangular trust region. 
The trust region is defined similar to $\Omega_k$ in the global sampling approach with two critical modifications. First, the region is based around $(x_k,y_k)$ instead of $(x_0,y_0)$, and second, we multiply the region by a confidence factor $\gamma\in[0,1]$ to allow for the size of the region to depend on the overall confidence in the new direction vector rather than the step size. Together these changes make up the region $\Omega_k=\{(x, y)\in(1-\gamma)[x_{k-1}-x_{min},x_{max}-x_{k-1}]\times(1-\gamma)[y_{k-1}-y_{min},y_{max}-y_{k-1}]: k=1...N\}$. As $\gamma\to1$, we are more confident in the updated direction so the search region for the next step can be reduced in size. Conversely, as $\gamma\to0$, we are less confident in the direction so the search size approaches the full parameter space. To prevent the full parameter space from being used on the first step, $\gamma$ is initially set close to 1. 
To update the confidence factor, we use the component-wise standard deviations of the direction vectors of the previous five steps. Because the direction vectors are unit vectors, the largest standard deviation is bounded at one in the case that 50\% of the points are at $-1$ and the other 50\% are at 1. 
For a general system with $D$-dimensional parameter space $\vec{\mu}=(\mu_1,...,\mu_D)$ the confidence factor can be computed as $\gamma = 1 - \prod_{i=1}^D\sigma_i^{(p)}$ where $\sigma_i$ is the standard deviation of the previous $p$ direction vectors of component $i$. Performing this updating algorithm on the Lorenz system from the same starting point, the path shown in Fig.~\ref{fig:lorenz_paths}~(b) is obtained where the blue region around the path shows the significantly smaller sampling region used by the navigation scheme to generate the path. 

We see from these preliminary results that using derivative free optimization works quite well for finding a periodic solution of this system, but it required sampling points in the vicinity of the current path point. Sampling the loss function is very expensive in this case because it requires a full numerical simulation of the system. This motivates the need for generating these paths using gradient descent to minimize the number of loss function samples required for moving through the space. The trade-off here is computing the gradient of the loss function. However, this is much more practical for driving a physical system because it only requires small changes in the parameters to estimate the gradient.

\section{Numerical Validation}\label{sec:param_opt_experiments}

In this section, we show numerical studies using three dynamical systems to demonstrate the capabilities and limitations of this method. Note that all results utilize the Adam optimizer for performing gradient descent to leverage the momentum advantages for avoiding local minima in the loss functions.

\subsection{Rössler System}

\begin{figure}[htbp]
  \centering
  \begin{minipage}[b]{0.32\textwidth}
    \centering
    \includegraphics[width=\textwidth]{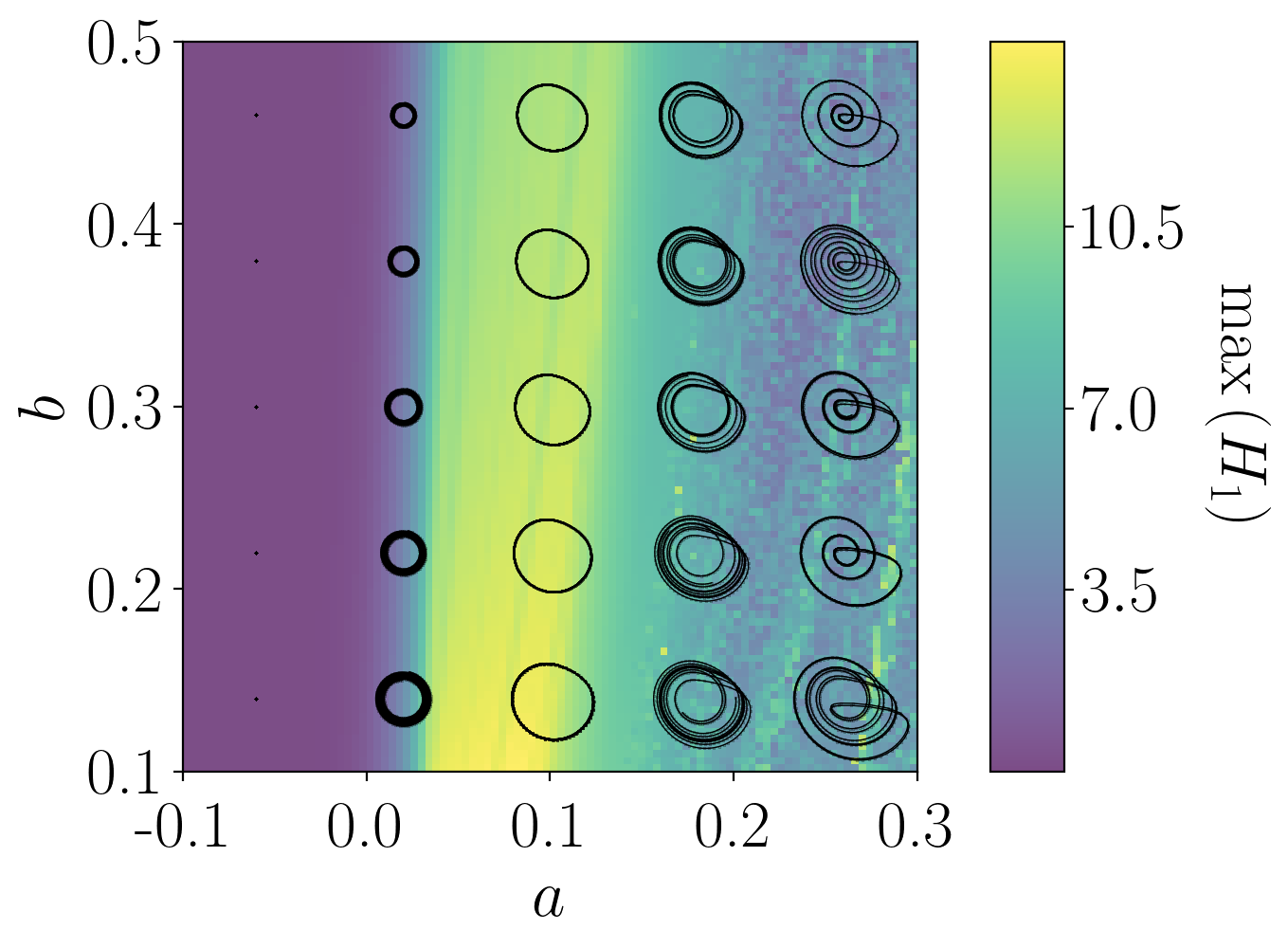}
    (a) Maximum Persistence
  \end{minipage}
  \begin{minipage}[b]{0.32\textwidth}
    \centering
    \includegraphics[width=\textwidth]{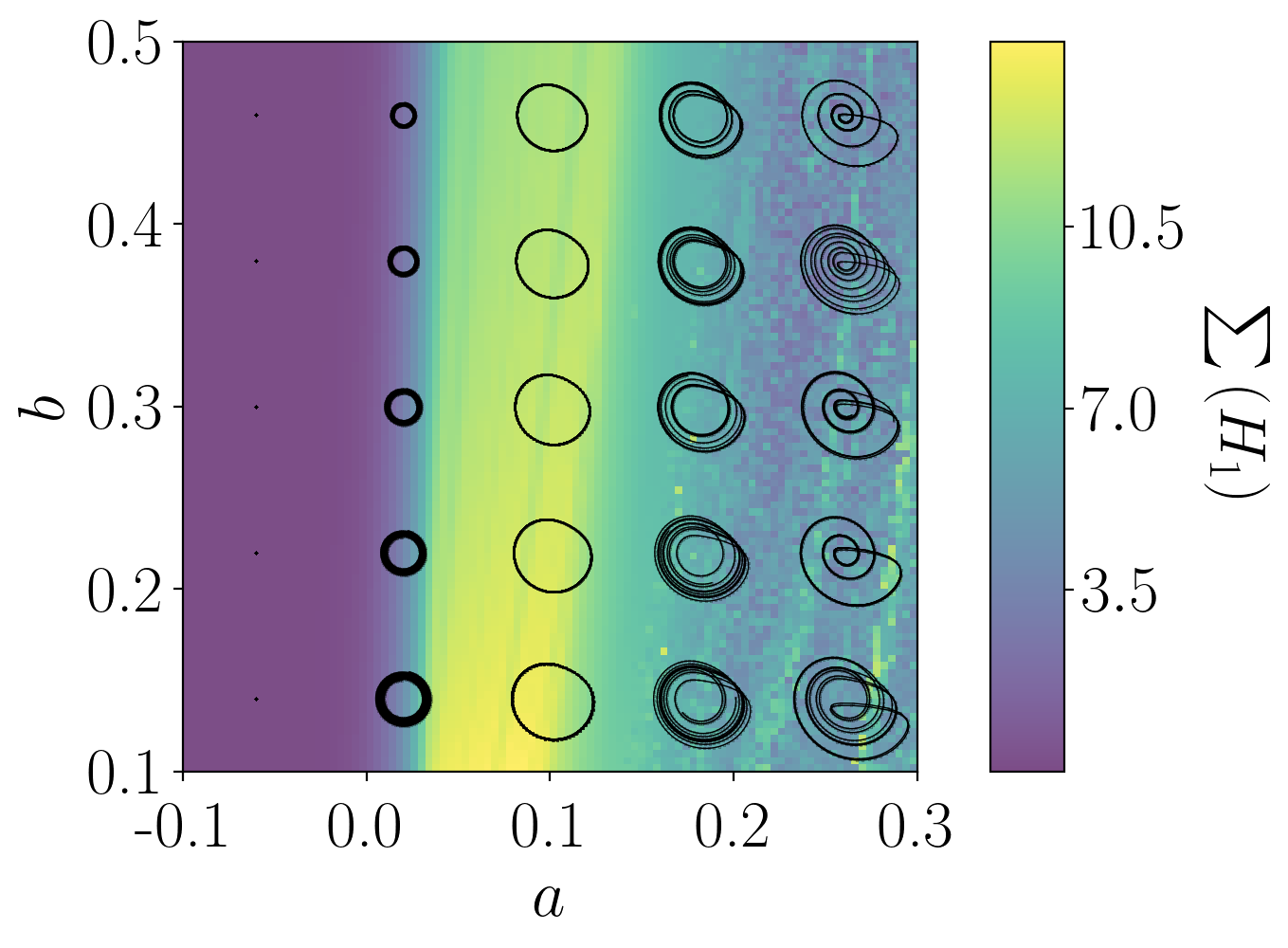}
    (b) Total Persistence
  \end{minipage}
  \begin{minipage}[b]{0.32\textwidth}
    \centering
    \includegraphics[width=\textwidth]{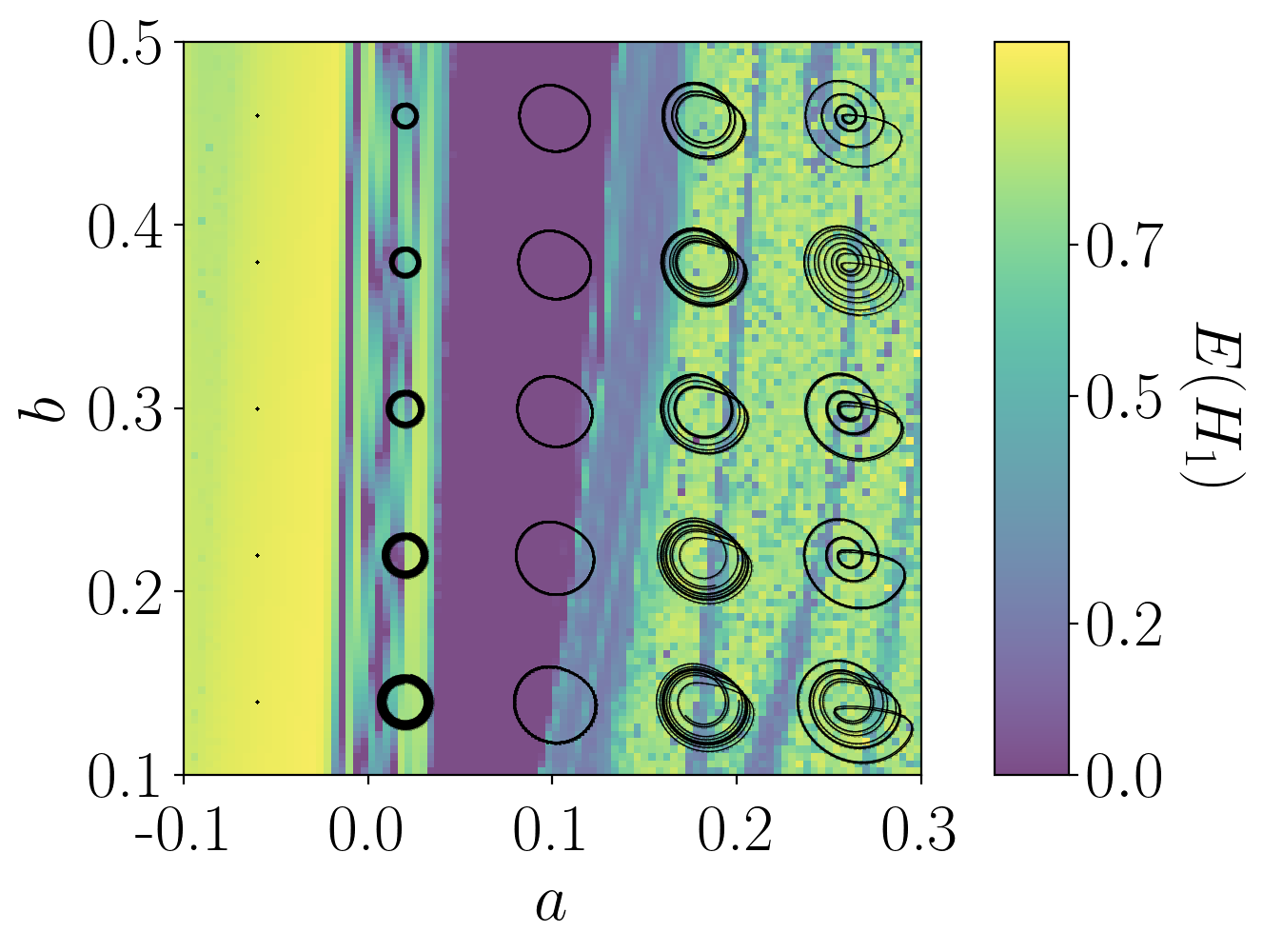}
    (c) Persistent Entropy
  \end{minipage}
  \caption{Rössler system persistence features plotted over a range of $a$ and $b$ values.}
  \label{fig:rossler_features}
\end{figure}

The first system we studied is the Rössler system, which is a well understood chaotic system \cite{Genesio_2008} described by the following set of differential equations:
\begin{equation}
\begin{aligned}
  \dot{x} &= -y - z, \\
  \dot{y} &= x + a y, \\
  \dot{z} &= b + z (x - c),
\end{aligned}
\end{equation}
where $a$, $b$, and $c$ are system parameters. For this analysis, we chose to vary $a$ and $b$ while keeping $c$ fixed at 5.7. The system was simulated using an initial condition of $(x, y, z) = (-0.4, 0.6, 1)$ and integrated over a time span from 0 to 200 time units, sampling every 0.04 time units and taking the last 500 points as the steady state response. Chaotic systems present a difficulty in that the trajectories deviate exponentially for an infinitesimal change in parameters. In theory this should be avoided because the overall topology of the trajectory should remain similar. However, since persistence pairs are being differentiated and not the overall shape or distribution of the pairs, the gradients explode and inconsistently vary by multiple orders of magnitude in chaotic regions of the parameter space. In practice, this issue is commonly alleviated with gradient clipping where the norm of the gradient is saturated at a specified magnitude and this has been shown to greatly improve training and exploding gradient issues in machine learning \cite{zhang2019}. The maximum persistence, total persistence and normalized persistent entropy are plotted over a range of $a$ and $b$ values to identify regions of chaotic and periodic behavior as shown in Fig.~\ref{fig:rossler_features}. We see that for larger $a$ values the system appears to be chaotic and as $a$ decreases it moves to periodic and fixed point responses. In the chaotic region of the parameter space, all three features appear to vary significantly suggesting that it will be difficult to find an optimal path in that region. This is further complicated by the entropy being high in the fixed point region. Even though the entropy was normalized, it is still at a maximum value in this region due to the significant number of low-lifetime persistence features. Intuitively, the entropy in this region is small because the loops are so small but in order for this to be detected with persistence it requires very long simulation time and makes computing the gradient using this pipeline computationally expensive. Nonetheless, we generated different results for this system to visualize and attempt to find optimal paths in this parameter space. 

\subsubsection{Chaos $\to$ Periodic}

\begin{figure}[htbp]
  \centering
  \includegraphics[width=\textwidth]{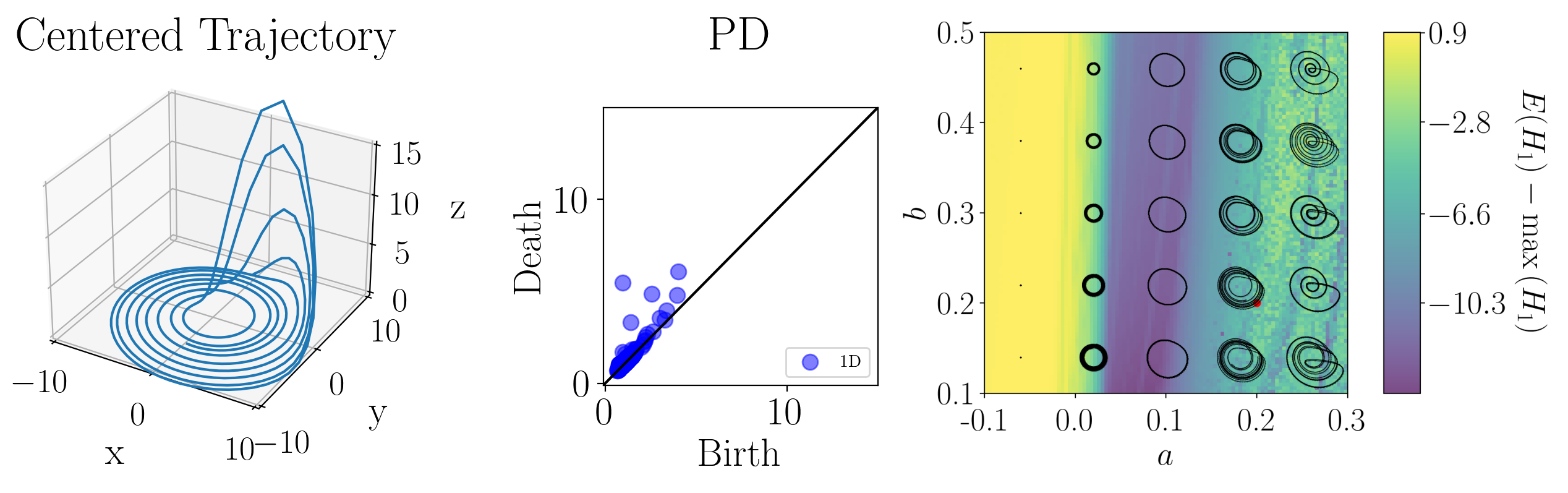}
  \caption{Initial Rössler trajectory starting at $a=0.2$ and $b=0.2$.}
  \label{fig:rossler_chaos_lc_initial_0.99}
\end{figure}

\begin{figure}[htbp]
  \centering
  \includegraphics[width=\textwidth]{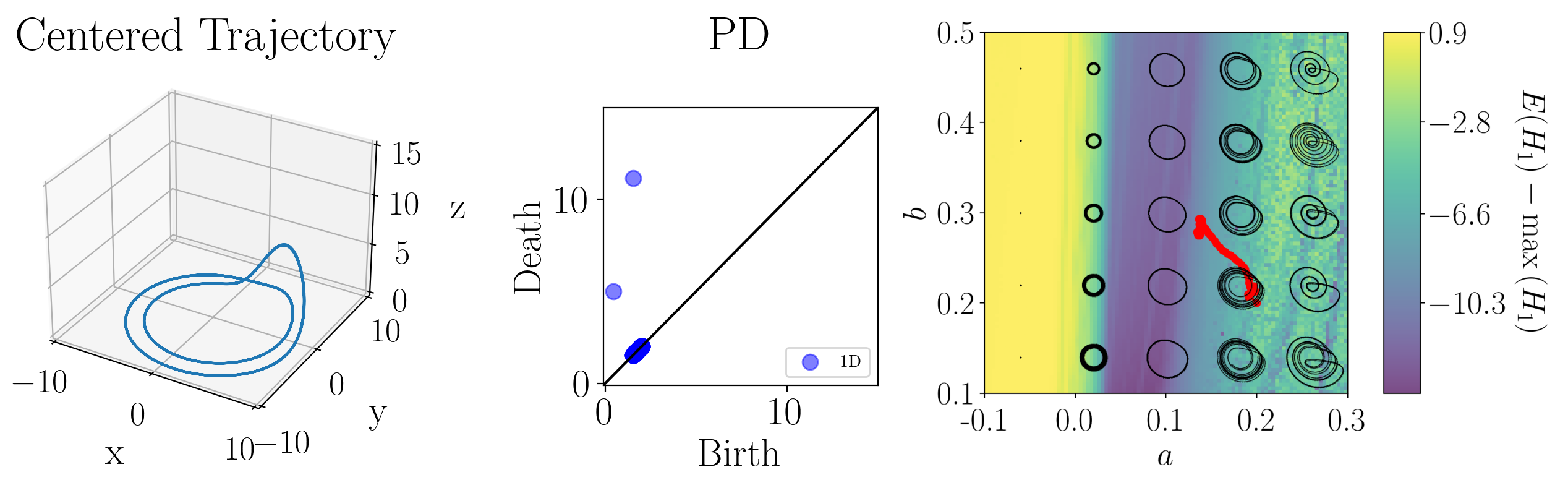}
  \caption{Final Rössler trajectory and parameter path after 225 epochs minimizing the persistent entropy and maximizing the maximum persistence using a learning rate decay of 0.99.}
  \label{fig:rossler_chaos_lc_final_0.99}
\end{figure}

For the first example, we attempted to move from the chaotic region of the parameter space to find a periodic solution using the persistent entropy and maximum persistence features. Specifically, we set the loss function to minimize persistent entropy and maximize maximum persistence ($\mathcal{L}=E-\text{maxPers}_1$). This leads to another complication from solving this multi-objective optimization problem. We see from Fig.~\ref{fig:rossler_features} that the maximum persistence and entropy are on different scales so optimizing this loss function will likely lead to an incorrect solution because it is not balanced. A common approach for avoiding this issue is to introduce scaling for each feature in the loss function as $\mathcal{L}=\sum_i \lambda_i \mathcal{L}_i$ where each $\lambda_i$ is a scale factor applied to the loss for each objective \cite{bischof2021}. For this work, we divide each persistence feature by its current value (detached from the gradient with pytorch) to ensure that each loss is on the same scale before being differentiated. This technique is presented in \cite{Xiao2024}. With this method, the value of each loss term will always be one, but the terms are all balanced and the gradients still vary in magnitude. We started the path at $a=0.2$, $b=0.2$ and used a learning rate of 0.01 with a gradient norm clip of 1.0. Because of the learning rate clip, it is critical to apply learning rate decay to make sure the solution converges so a decay of 1\% per epoch was applied. The initial and final results are shown in Figs.~\ref{fig:rossler_chaos_lc_initial_0.99} and \ref{fig:rossler_chaos_lc_final_0.99} where we see that the path successfully exited the chaotic region and approached a periodic solution, but due to the learning rate decay the optimizer does not explore the parameter space enough and it settles on a more complicated periodic solution with two prominent loops in the persistence diagram. 

\begin{figure}[htbp]
  \centering
  \includegraphics[width=\textwidth]{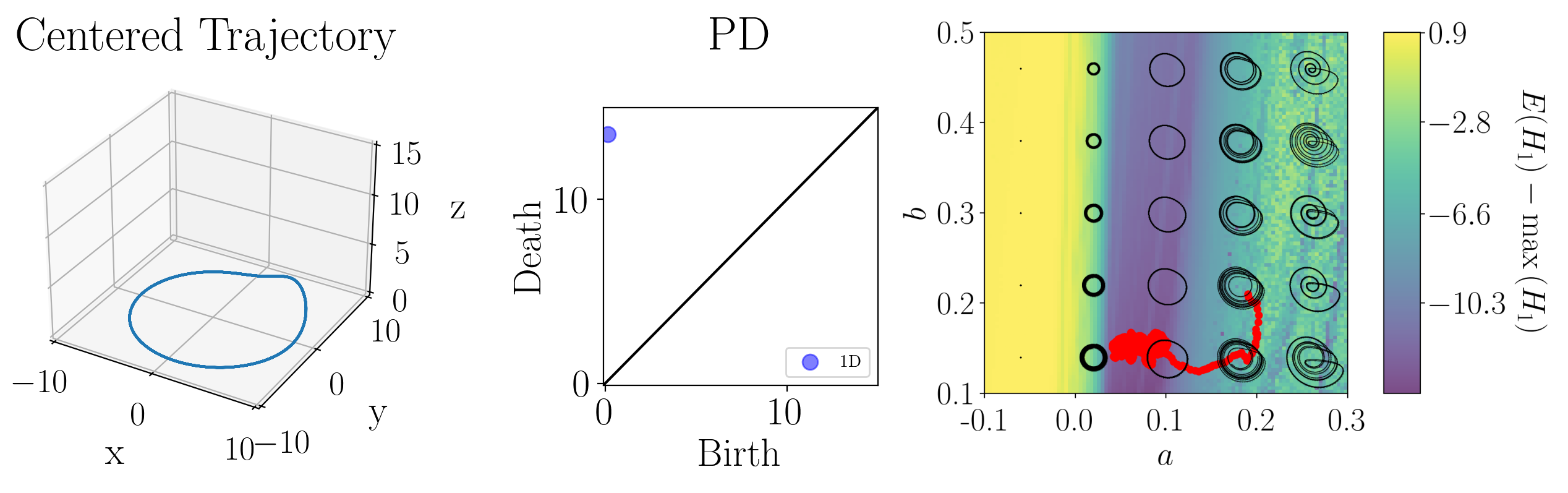}
  \caption{Final Rössler trajectory and parameter path after 450 epochs minimizing the persistent entropy and maximizing the maximum persistence using a learning rate decay of 0.999.}
  \label{fig:rossler_chaos_lc_final_0.999}
\end{figure}

For the next test, we used an identical setup, but change the learning rate decay rate to 0.999 to allow for more exploration and the resulting path is shown in Figs.~\ref{fig:rossler_chaos_lc_final_0.999}. We see that the path converged on a simpler periodic solution in this case, but required many more steps due to the slower learning rate decay. Interestingly, the path seemed to oscillate in the periodic region of the parameter space without entering the fixed point region or chaotic region again. Due to starting in the chaotic region the resulting path is also drastically different from the other learning rate decay which demonstrates the difficulty of optimizing when responses are chaotic.

\subsubsection{Periodic $\to$ Fixed Point}

\begin{figure}[htbp]
  \centering
  \includegraphics[width=\textwidth]{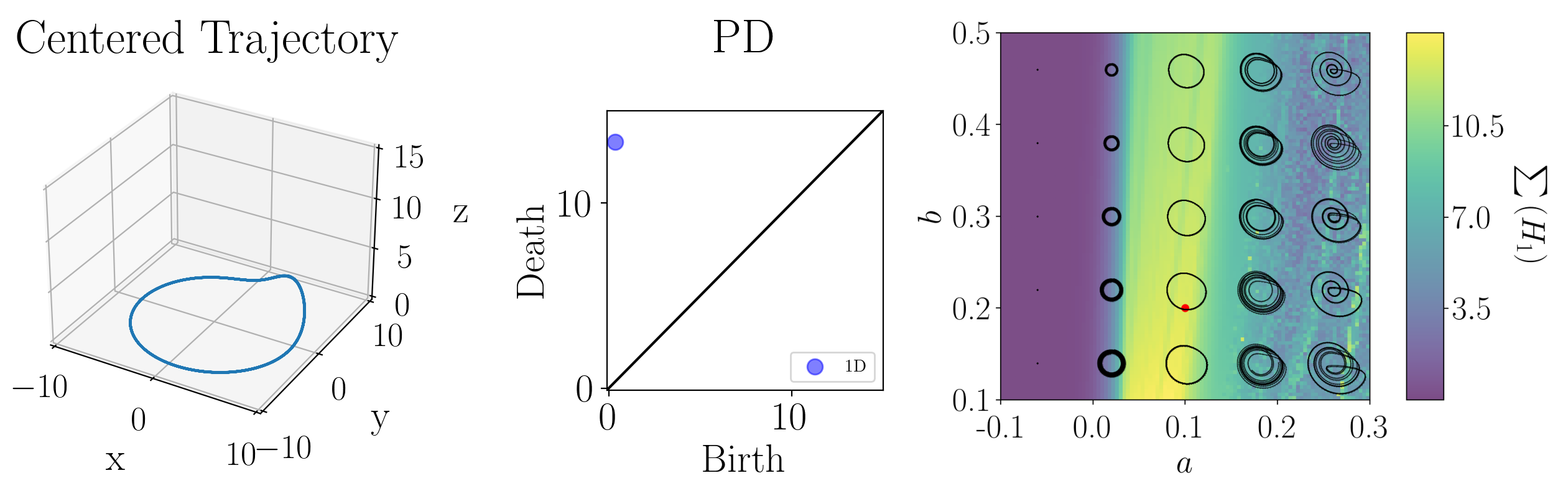}
  \caption{Initial Rössler trajectory at $a=0.1$ and $b=0.2$.}
  \label{fig:rossler_LC_FP_initial0.01_0.99}
\end{figure}

\begin{figure}[htbp]
  \centering
  \includegraphics[width=\textwidth]{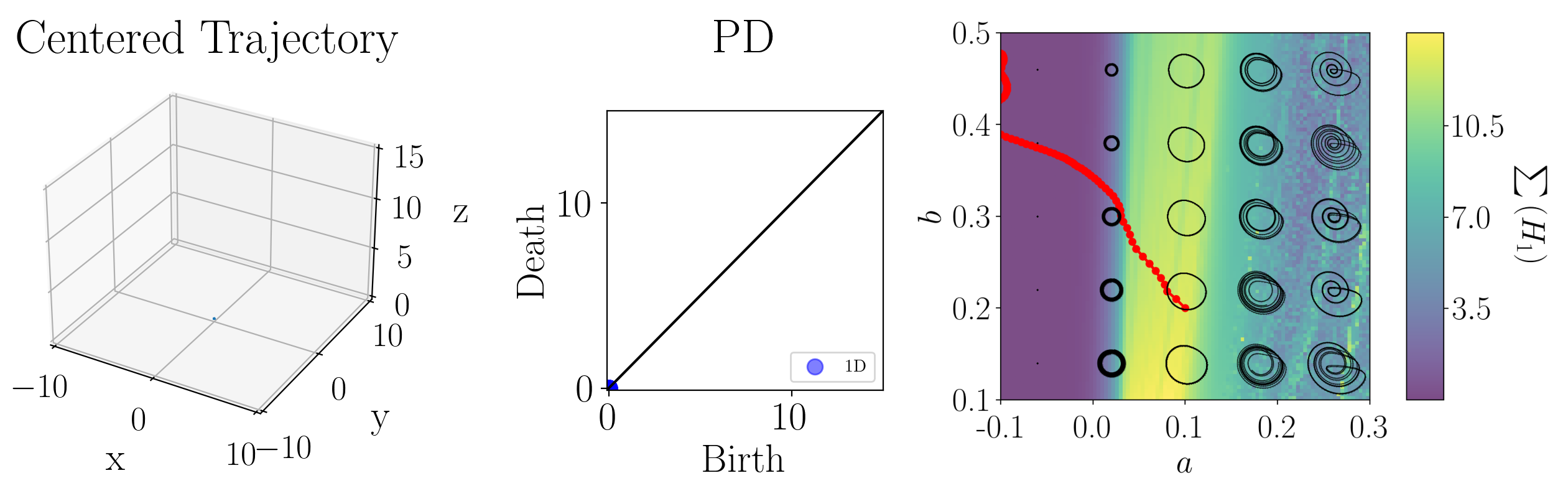}
  \caption{Final Rössler trajectory and parameter path after 142 epochs minimizing the total persistence using a learning rate decay of 0.99.}
  \label{fig:rossler_LC_FP_final0.01_0.99}
\end{figure}

Next we show examples attempting to move from a periodic response to fixed point response using the features from Section~\ref{sec:cost_lib}. Initially, the goal was to minimize the total persistence to reach the desired response. We set the cost function to be the total persistence and started the path at $a=0.1$ and $b=0.2$. The initial trajectory and persistence diagram are shown in Fig.~\ref{fig:rossler_LC_FP_initial0.01_0.99}. Using the total persistence cost function and a learning rate of 0.01 with decay rate of 0.99, the path in Fig.~\ref{fig:rossler_LC_FP_final0.01_0.99} was obtained after 142 epochs. We see that the path correctly moves toward the fixed point region but continues to hit the lower bound on $a$. The path oscillates on this line before converging on parameters slightly outside of the region. This result demonstrates a limitation of this method in the way the Adam optimizer works. Because the optimizer had momentum moving toward the wall the path was able to exit the region. This is why when defining forbidden regions using the method in Section~\ref{sec:forbidden}, it is important to give some buffer space as it is possible for the path to go beyond the boundary in some cases.

\begin{figure}[htbp]
  \centering
  \includegraphics[width=\textwidth]{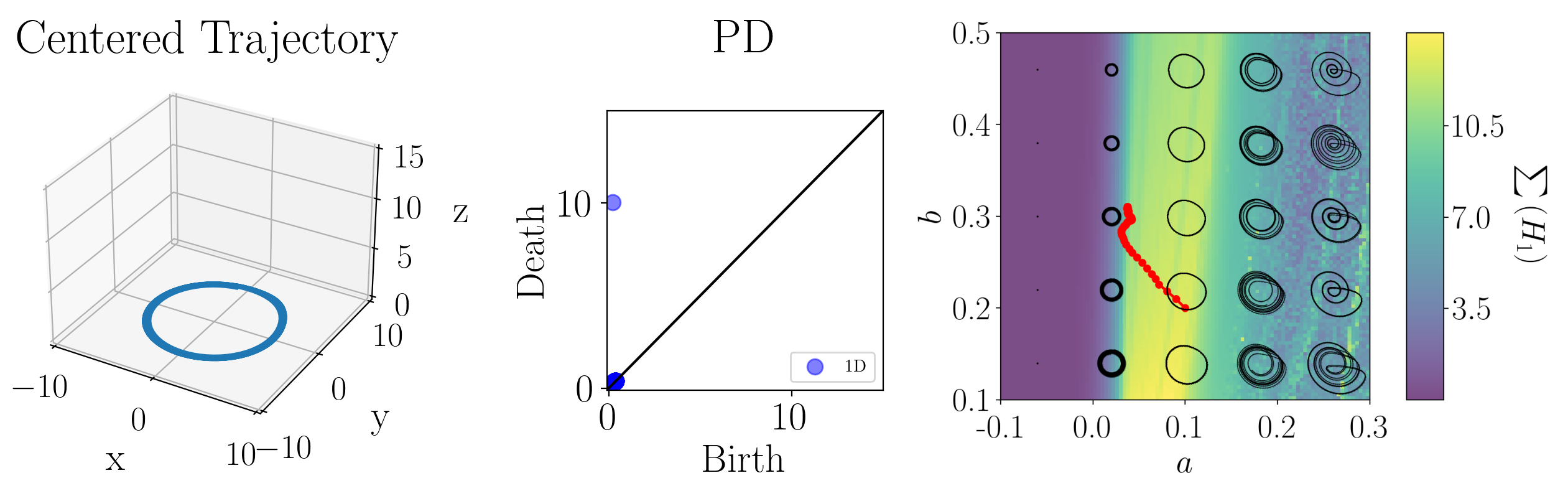}
  \caption{Final Rössler trajectory and parameter path after 106 epochs minimizing the total persistence using a learning rate decay of 0.95.}
  \label{fig:rossler_LC_FP_final0.01_0.95}
\end{figure}

\begin{figure}[htbp]
  \centering
  \includegraphics[width=\textwidth]{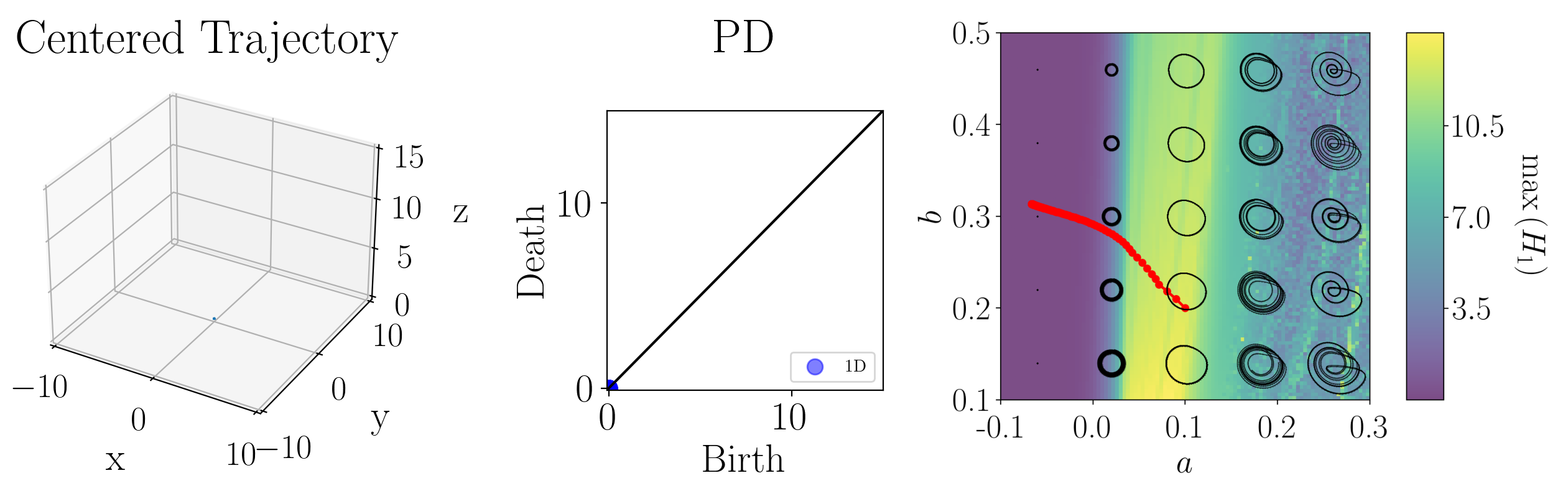}
  \caption{Final Rössler trajectory and parameter path after 147 epochs minimizing the maximum persistence using a learning rate decay of 0.95.}
  \label{fig:rossler_LC_FP_final0.01_0.95_MP}
\end{figure}

Next, we ran the same test but reduced the learning rate decay to 0.95 and the path in Fig.~\ref{fig:rossler_LC_FP_final0.01_0.95} was obtained after 106 epochs. We see that the total persistence loss landscape provides some resistance to the fixed point region because the total persistence slightly increases as more low lifetime persistence pairs appear. The path was never able to enter the fixed point region in this example. This demonstrates a limitation of the total persistence feature that is mentioned in Section~\ref{sec:cost_lib} where the total persistence can be biased by the number of pairs in the persistence diagram. To attempt to mitigate this issue, we chose to run the same test using maximum persistence and the results are shown in Fig.~\ref{fig:rossler_LC_FP_final0.01_0.95_MP} where we see the path converges on the parameters $a=-0.0664$ and $b=0.3136$ after 147 epochs.

\subsubsection{Chaos $\to$ Fixed Point}

\begin{figure}[htbp]
  \centering
  \includegraphics[width=\textwidth]{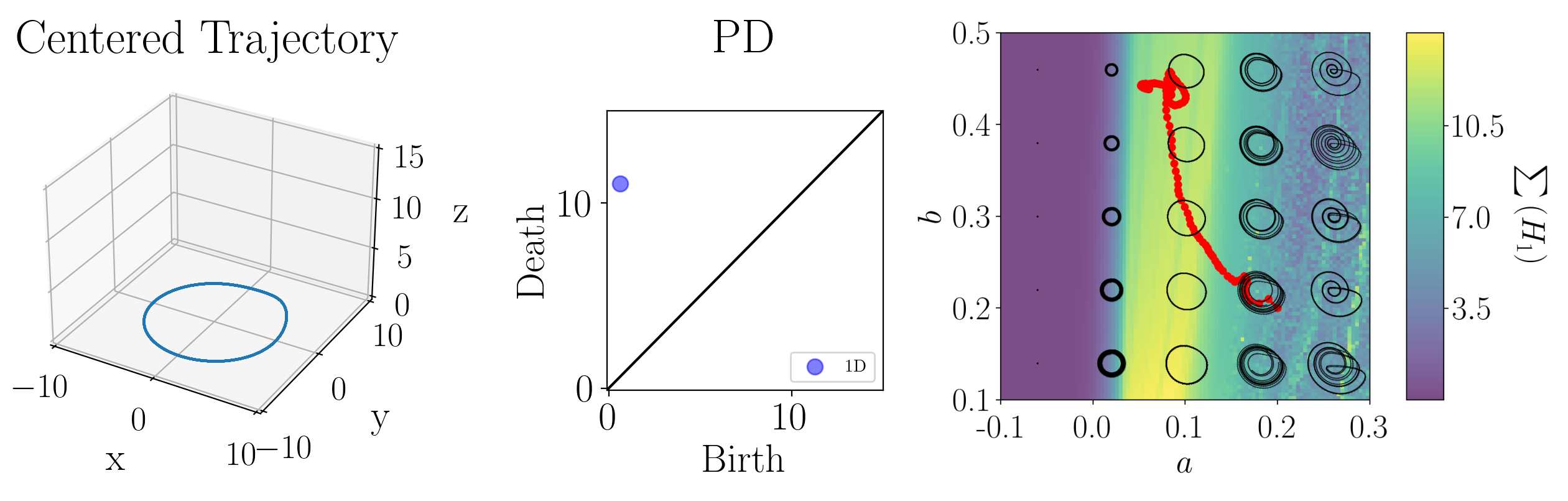}
  \caption{Final Rössler trajectory and parameter path after 165 epochs minimizing the maximum persistence using a learning rate of 0.01 and decay rate of 0.999.}
  \label{fig:rossler_chaos_fp_initial_0.999}
\end{figure}

\begin{figure}[htbp]
  \centering
  \includegraphics[width=\textwidth]{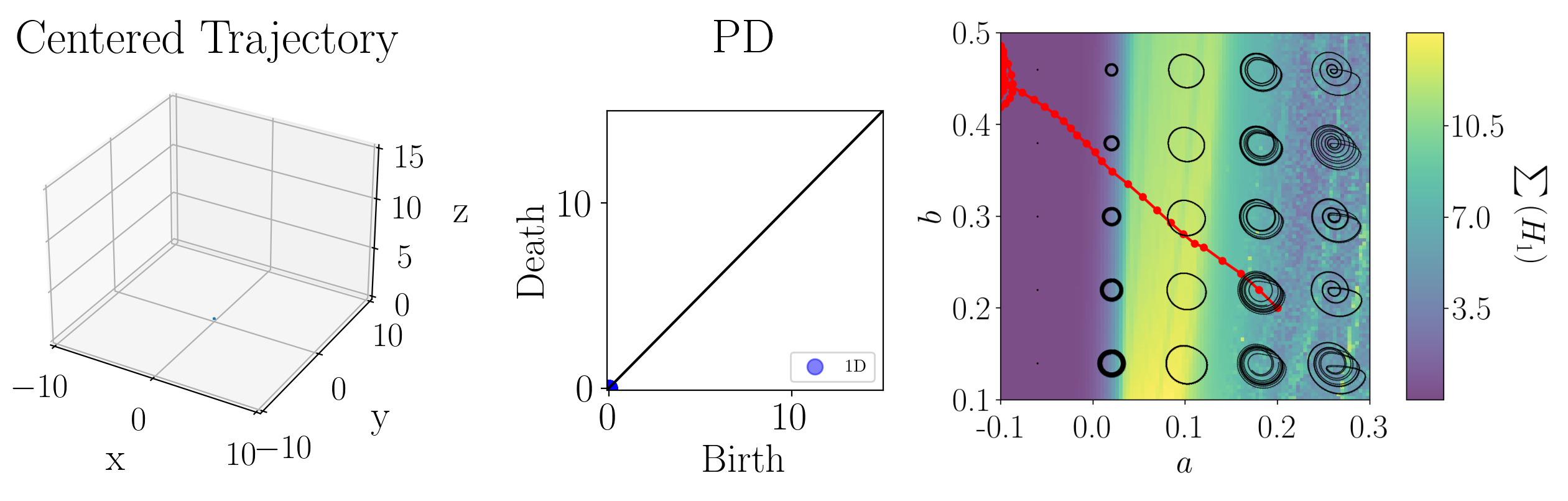}
  \caption{Final Rössler trajectory and parameter path after 200 epochs minimizing the maximum persistence using a learning rate of 0.02 and decay rate of 0.999.}
  \label{fig:rossler_chaos_fp_initial_0.999_0.02}
\end{figure}

For the final Rössler system test, we aimed to move from the chaotic region to a fixed point response. We used the same setup as the previous tests by starting the path at $a=0.2$, $b=0.2$ and set the learning rate to 0.01 with a decay rate of 0.999. Because the persistent entropy reaches a minimum in the periodic region, we chose to only use total persistence in this case for the loss function. The initial trajectory is shown in Fig.~\ref{fig:rossler_chaos_lc_initial_0.99} and after 165 epochs of optimization the path is shown in Fig.~\ref{fig:rossler_chaos_fp_initial_0.999}. The path in this case correctly exited the chaotic region of the parameter space, but was not able to enter the fixed point region likely due to momentum issues again. Once the path entered the periodic region, interestingly, it moved in the vertical $b$ direction and it is believed that this is also a consequence of the optimizers momentum. The gradient in the $a$ direction is much larger than the gradient in the $b$ direction, so the Adam optimizer conservatively moves vertically and due to the small learning rate it is not able to enter the fixed point region. To fix this issue, we increased the learning rate to 0.02 and ran the same test. The resulting path is shown in Fig.~\ref{fig:rossler_chaos_fp_initial_0.999_0.02}. The resulting path successfully reached a fixed point solution from chaos with the larger learning rate, but also required regularization as the optimizer tried to decrease $a$ beyond the lower limit.

\subsection{Magnetic Pendulum}

\begin{figure}[htbp]
  \centering
  \includegraphics[width=0.5\textwidth]{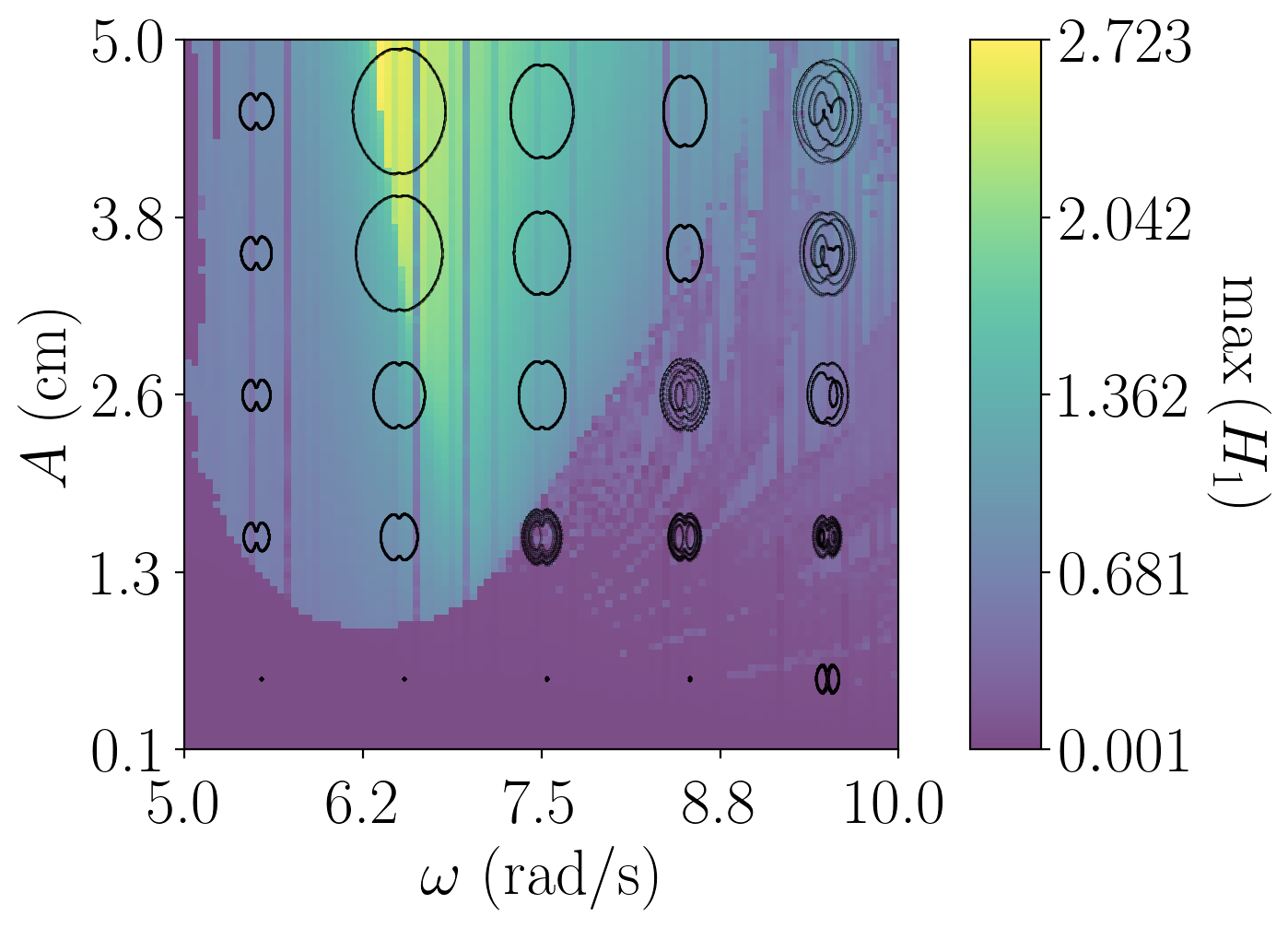}
  \caption{Maximum persistence plotted over a range of base excitation amplitude and frequency for the magnetic pendulum system.}
  \label{fig:mag_pen_maxpers}
\end{figure}

\begin{figure}[htbp]
  \centering
  \includegraphics[width=\textwidth]{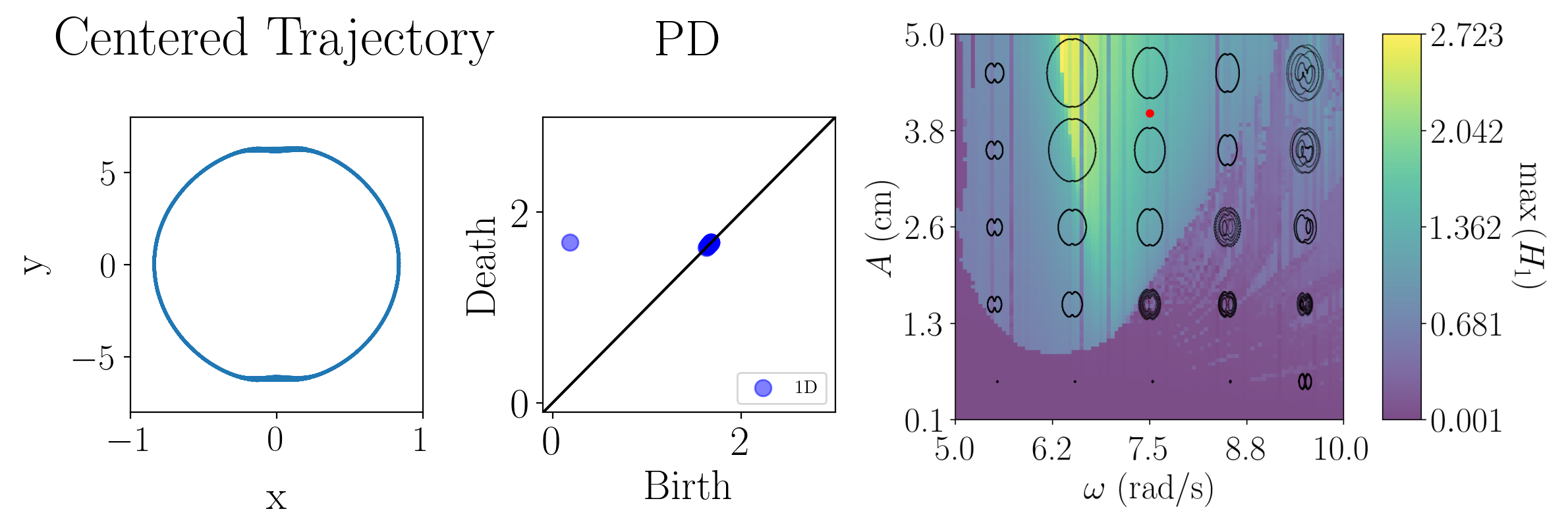}
  \caption{Initial response and persistence diagram for the magnetic pendulum system.}
  \label{fig:mag_pen_initial}
\end{figure}

\begin{figure}[htbp]
  \centering
  \includegraphics[width=\textwidth]{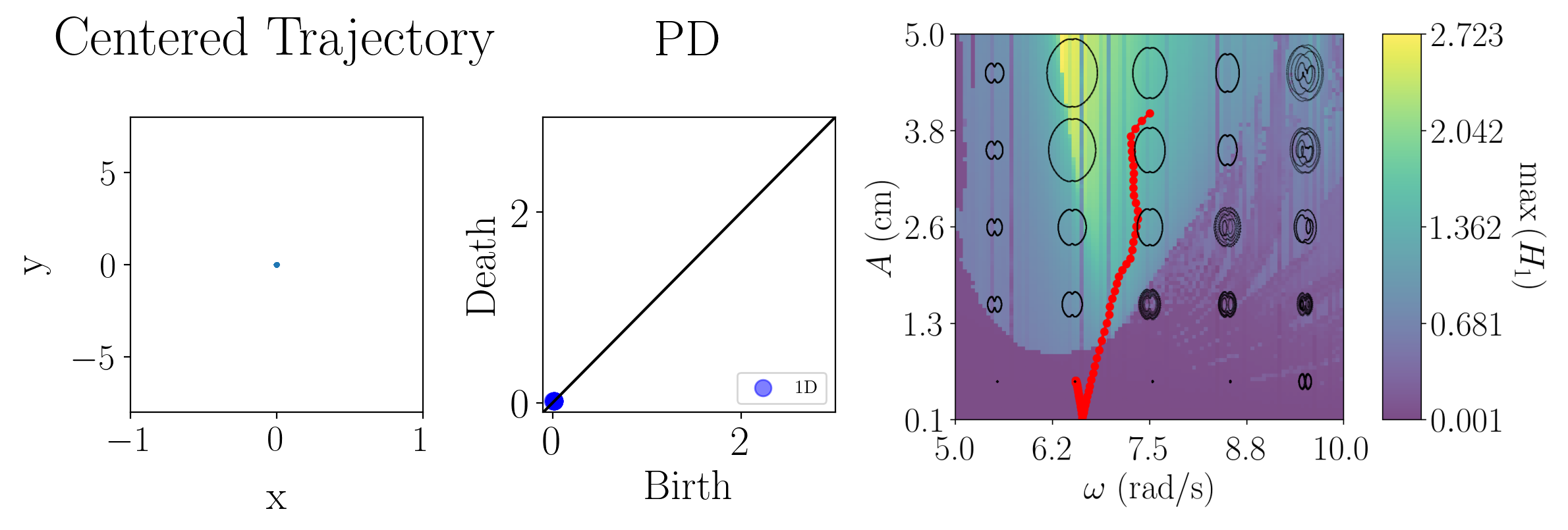}
  \caption{Final response and persistence diagram for the magnetic pendulum system with the parameter space path minimizing maximum persistence.}
  \label{fig:mag_pen_final}
\end{figure}

Next, we studied the base excited magnetic pendulum system shown in \cite{Myers_2020_pendulum}. This system consists of a normal base excited pendulum system with a magnet on the end of the pendulum arm and another magnet on the base. This creates a magnetic interaction between the pendulum and base that leads to complex behavior. The equation of motion for the pendulum is:
\begin{equation}
    (M r_{\text{cm}}^2 + I_{\text{cm}})\ddot{\theta}+M g r_{\text{cm}}\sin{(\theta)}  = - \tau_v - \tau_m - M r_{\text{cm}} \ddot{x}_{\text{base}} \cos\theta.
\end{equation}
where $\theta$ is the pendulum angle measured from vertical, $\ddot{x}_{\text{base}} = -A \omega^2 \sin(\omega t)$ is the base acceleration, $M = 0.1038$ kg is the total mass, $l = 0.208$ m is the length of the pendulum, $g = 9.81$ m/s$^2$, $r_{\text{cm}} = 0.18775$ m is the distance to the center of mass from the hinge, $I_{\text{cm}} = 1.919 \times 10^{-5}$ kg$\cdot$m$^2$ is the mass moment of inertia, $\mu_v = 0.003$ is the viscous damping coefficient, $m = 1.2$ A$\cdot$m$^2$ is the magnetic dipole moment, $\mu_0 = 1.257 \times 10^{-6}$ N/A$^2$ is the permittivity of free space, and $d = 0.032$ m is the minimum distance between the magnets. In \cite{Myers_2020_pendulum}, the authors measure system parameters and we use the same parameters in this work. $\tau_m$ is the magnetic interaction torque given by,
\begin{equation}
  \tau_m = l (F_r \cos(\phi - \theta) - F_\phi \sin(\phi - \theta)),
\end{equation}
where,
\begin{align}
  F_r &= \frac{3 \mu_0 m^2}{4 \pi r^4} \left( 2 \cos(\phi - a) \cos(\phi - b) - \sin(\phi - a) \sin(\phi - b) \right), \\
  F_\phi &= \frac{3 \mu_0 m^2}{4 \pi r^4} \sin(2\phi - a - b),
\end{align}
are the radial and tangential magnetic interaction forces, and $r$ and $\phi$ are the polar coordinate locations of the pendulum measured from the base magnet and they are computed with,
\begin{align}
  r &= \sqrt{l^2 + (d + l)^2 - 2l(l + d) \cos\theta}, \\
  \phi &= \frac{\pi}{2} - \arcsin\left(\frac{l}{r} \sin\theta\right),
\end{align}
with $a = \frac{3\pi}{2}$, $b = \frac{\pi}{2} - \theta$. Lastly, the viscous damping torque is,
\begin{align}
  \tau_v &= \mu_v \dot{\theta}.
\end{align}
For this analysis, we chose to vary the base excitation amplitude ($A$) and frequency ($\omega$) for persistence optimization. We plotted the maximum persistence over a range of amplitudes and frequencies shown in Fig.~\ref{fig:mag_pen_maxpers}. We see that for a range of larger amplitude and lower frequency, the response is periodic and for low amplitude and frequency the system approaches a fixed point. This intuitively makes sense for the system. The goal is to reach a fixed point using persistence optimization now. We set the cost function to minimize the maximum persistence and started the path at $A=4$ cm and $\omega=7.5$ rad/s. The system was simulated for 100 seconds sampling every 0.03 seconds and taking the last 500 points for computing persistence. The initial response and persistence diagram are shown in Fig.~\ref{fig:mag_pen_initial}. In this case, we wanted to explore a larger region of the parameter space, so we set the learning rate to 0.1. After 85 optimization steps, the path is shown in Fig.~\ref{fig:mag_pen_final}. We see that the path reaches a set of parameters that result in fixed point stability for the system and the regularization term caused the path to reflect off of the lower limit on $A$. Intuitively, we know that taking $A$ to be as small as possible is the best way to minimize oscillations, but in this case it settles on a set of somewhat nontrivial parameters at $A=0.57$ cm and $\omega=6.55$ rad/s. 


\subsection{Lorenz System}

\begin{figure}[htbp]
  \centering
  \begin{minipage}[b]{0.4\textwidth}
    \centering
    \includegraphics[width=\textwidth]{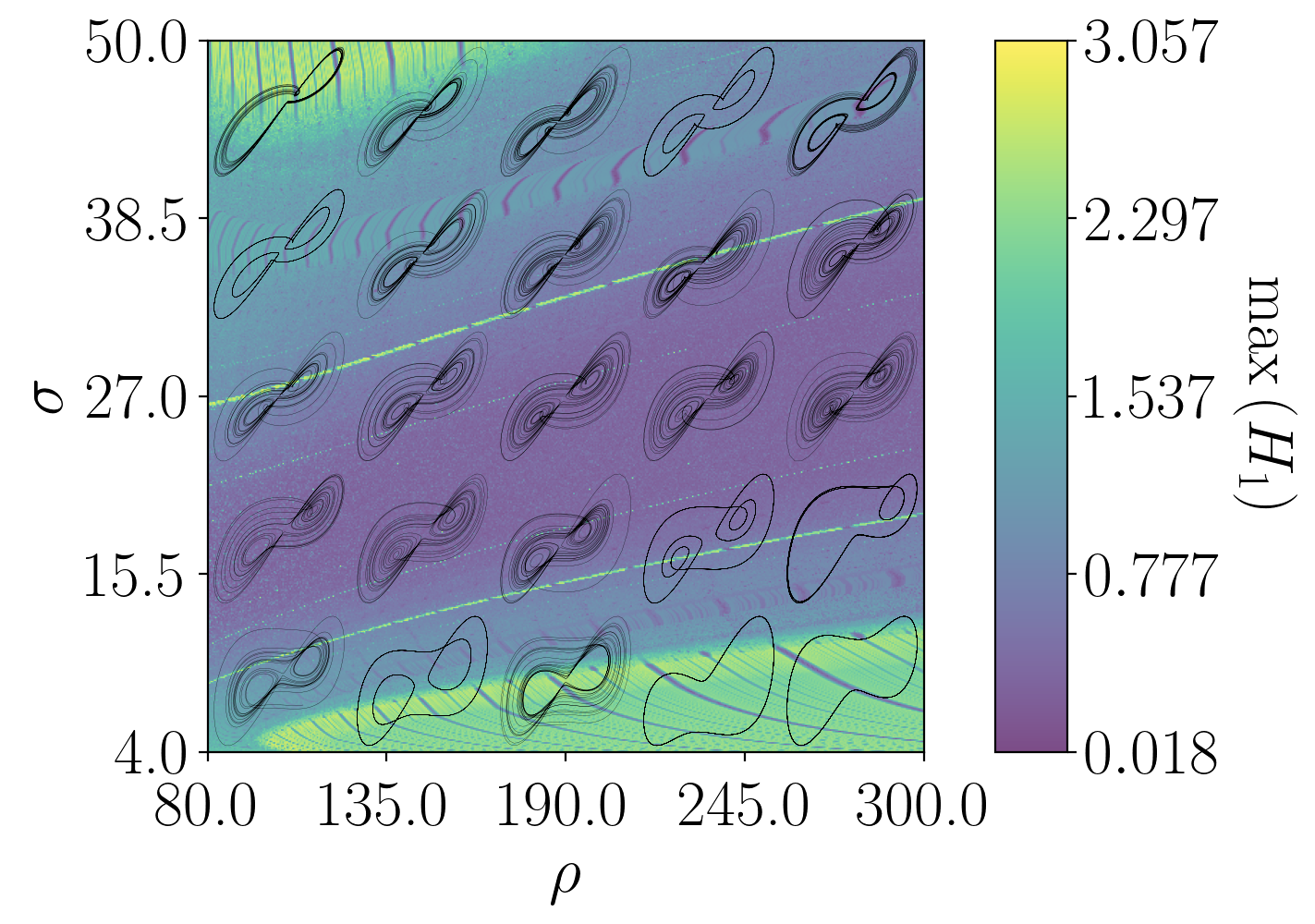}
    (a) Maximum Persistence
  \end{minipage}
  \begin{minipage}[b]{0.4\textwidth}
    \centering
    \includegraphics[width=\textwidth]{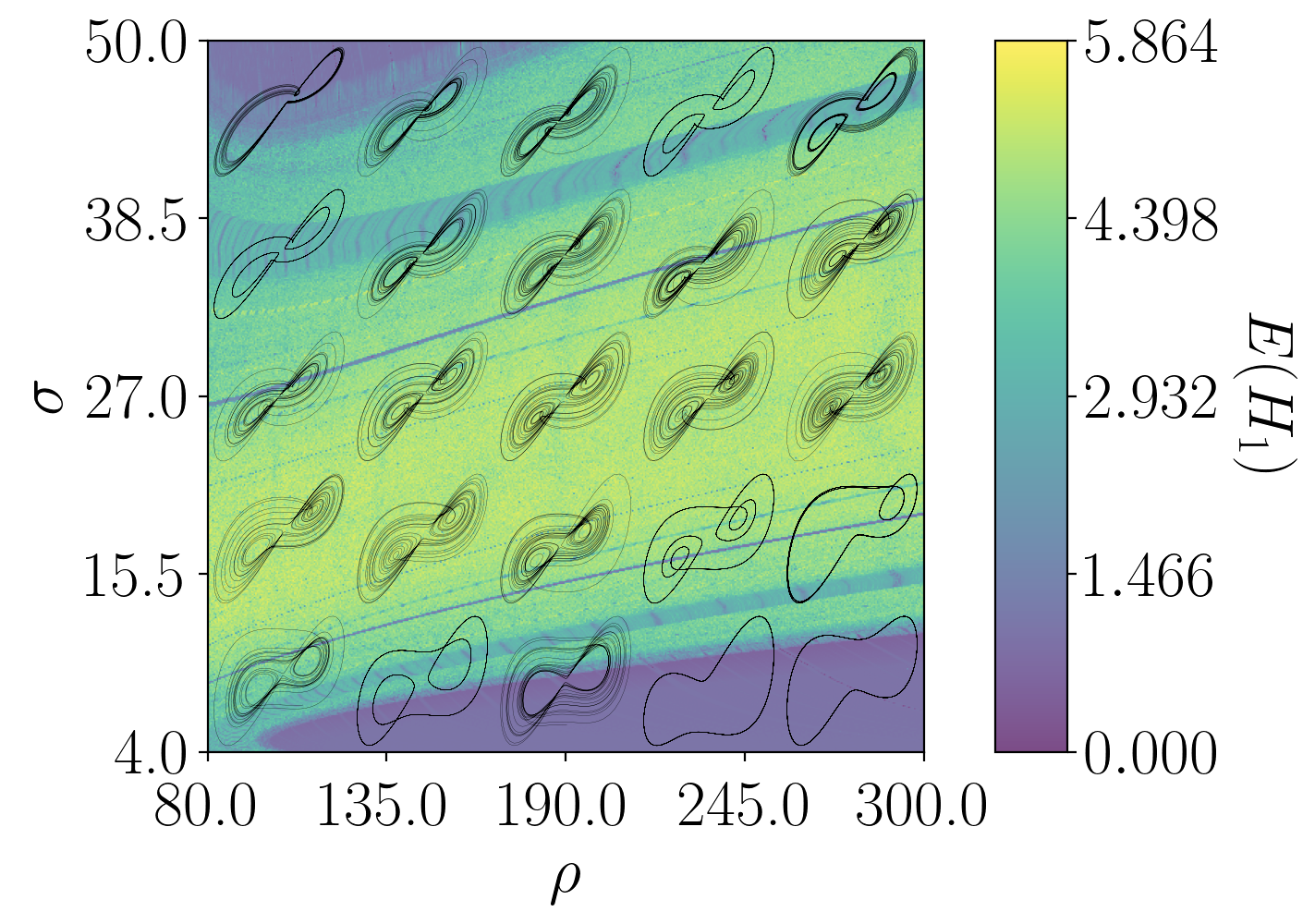}
    (b) Persistent Entropy
  \end{minipage}
  \caption{Lorenz system persistence features plotted over a range of $\rho$ and $\sigma$ values.}
  \label{fig:lorenz_features}
\end{figure}

\begin{figure}[htbp]
  \centering
  \includegraphics[width=\textwidth]{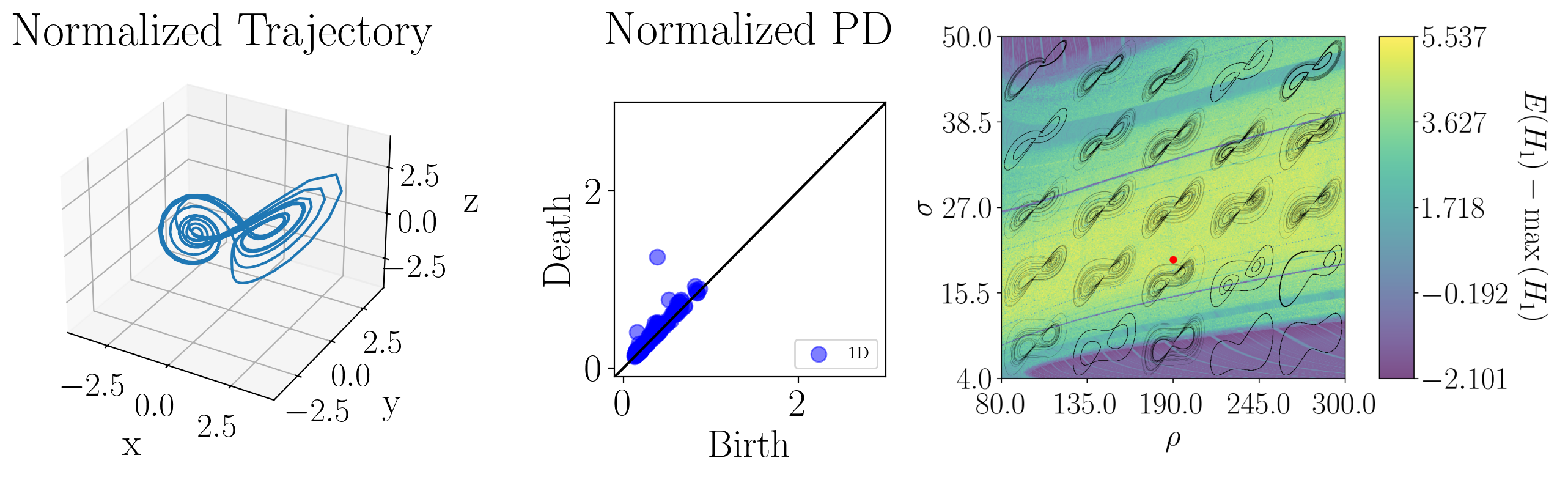}
  \caption{Initial Lorenz trajectory.}
  \label{fig:lorenz_initial}
\end{figure}

For the final example, returned to the Lorenz system to show that this method aligns with the preliminary results from Section~\ref{sec:preliminary_path}. For this analysis, we chose to vary $\sigma$ and $\rho$ while keeping $\beta$ fixed at the typical value of $8/3$ for visualization purposes. The system was simulated using an initial condition of $(x, y, z) = (1, 1, 1)$ and integrated over a time span from 0 to 10 time units, sampling every 0.01 time units and taking the last 500 points as the steady state response. The maximum persistence and persistent entropy were plotted over a range of $\sigma$ and $\rho$ values to identify regions of chaotic and periodic behavior as shown in Fig.~\ref{fig:lorenz_features}. We see that periodic solutions appear to be most prominent for low values of $\sigma$ and as this parameter increases, a chaotic region forms. The goal was to start the path in the chaotic region, and use these persistence features to navigate to the periodic region. However, this example is much more challenging than previous examples due to the parameter space being significantly larger. In all cases, the loss function was defined to be the difference between entropy and maximum persistence to minimize entropy and maximize maximum persistence and the loss function scaling was applied to ensure these terms were on the same scale. Gradient clipping to a norm of one was also applied in this example due to the exploding gradients in chaotic regions of the parameter space. We started by initializing the path at $\rho=190$ and $\sigma=20$. The initial trajectory is shown in Fig.~\ref{fig:lorenz_initial}. Note that for the plotting purposes, the trajectories were normalized using the mean and standard deviation, but for the persistence optimization computations the unmodified state space was used.  

\begin{figure}[htbp]
  \centering
  \includegraphics[width=\textwidth]{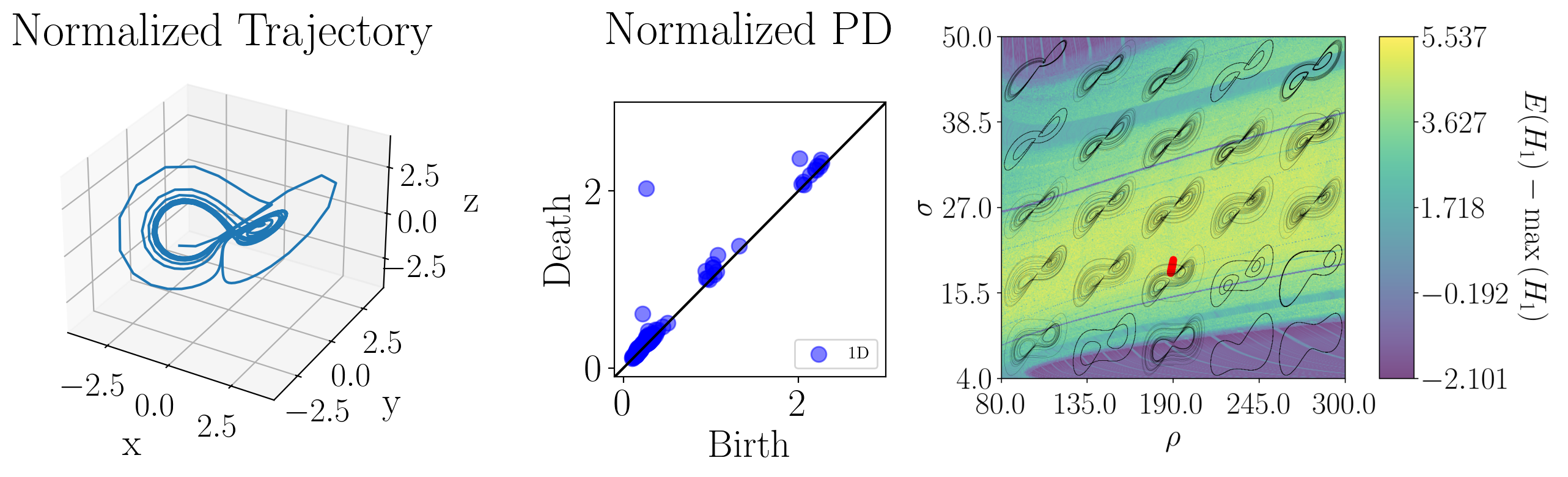}
  \caption{Final Lorenz trajectory after 385 epochs with a learning rate of 0.1 and no decay. }
  \label{fig:lorenz0.1_final}
\end{figure}

\begin{figure}[htbp]
  \centering
  \includegraphics[width=\textwidth]{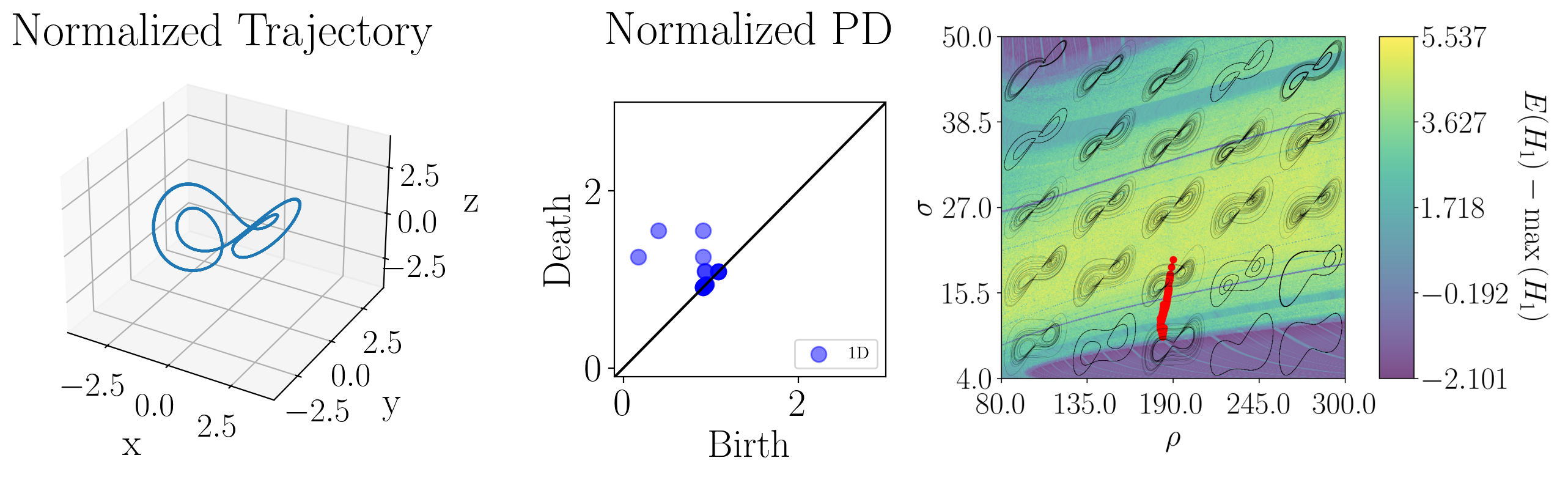}
  \caption{Final Lorenz trajectory after 382 epochs with a learning rate of 1.0 and a decay rate of 0.99. }
  \label{fig:lorenz1.0_0.99_final}
\end{figure}

\begin{figure}[htbp]
  \centering
  \includegraphics[width=\textwidth]{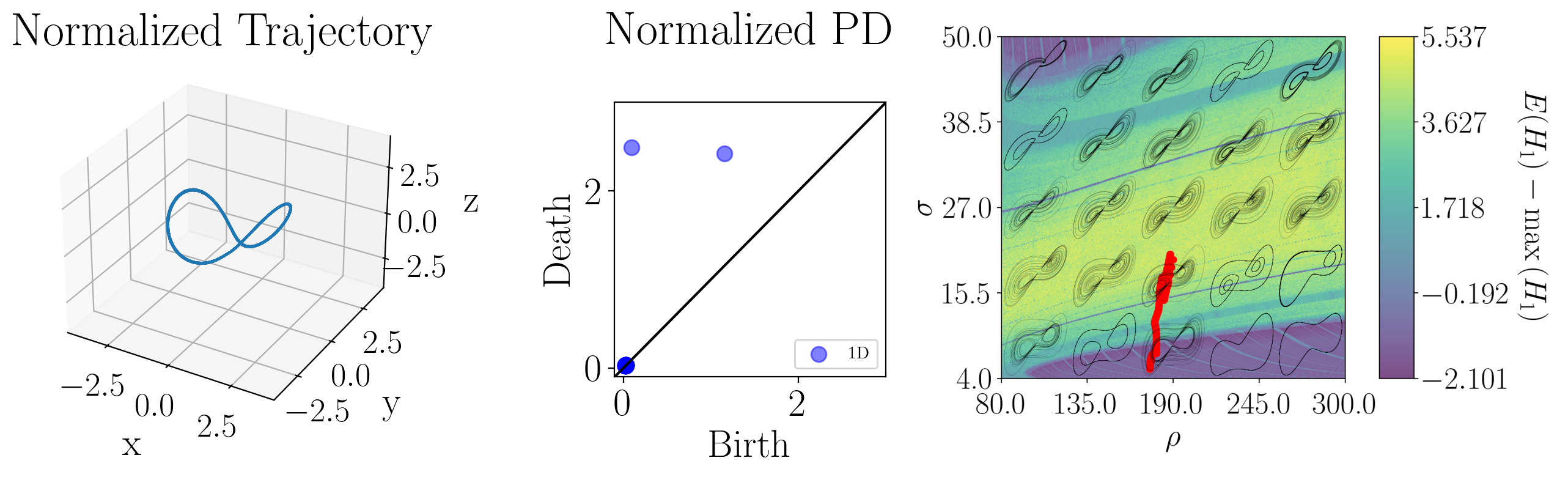}
  \caption{Final Lorenz trajectory after 2100 epochs with a learning rate of 1.0 and a decay rate of 0.995. }
  \label{fig:lorenz1.0_0.995_final}
\end{figure}

\begin{figure}[htbp]
  \centering
  \includegraphics[width=0.55\textwidth]{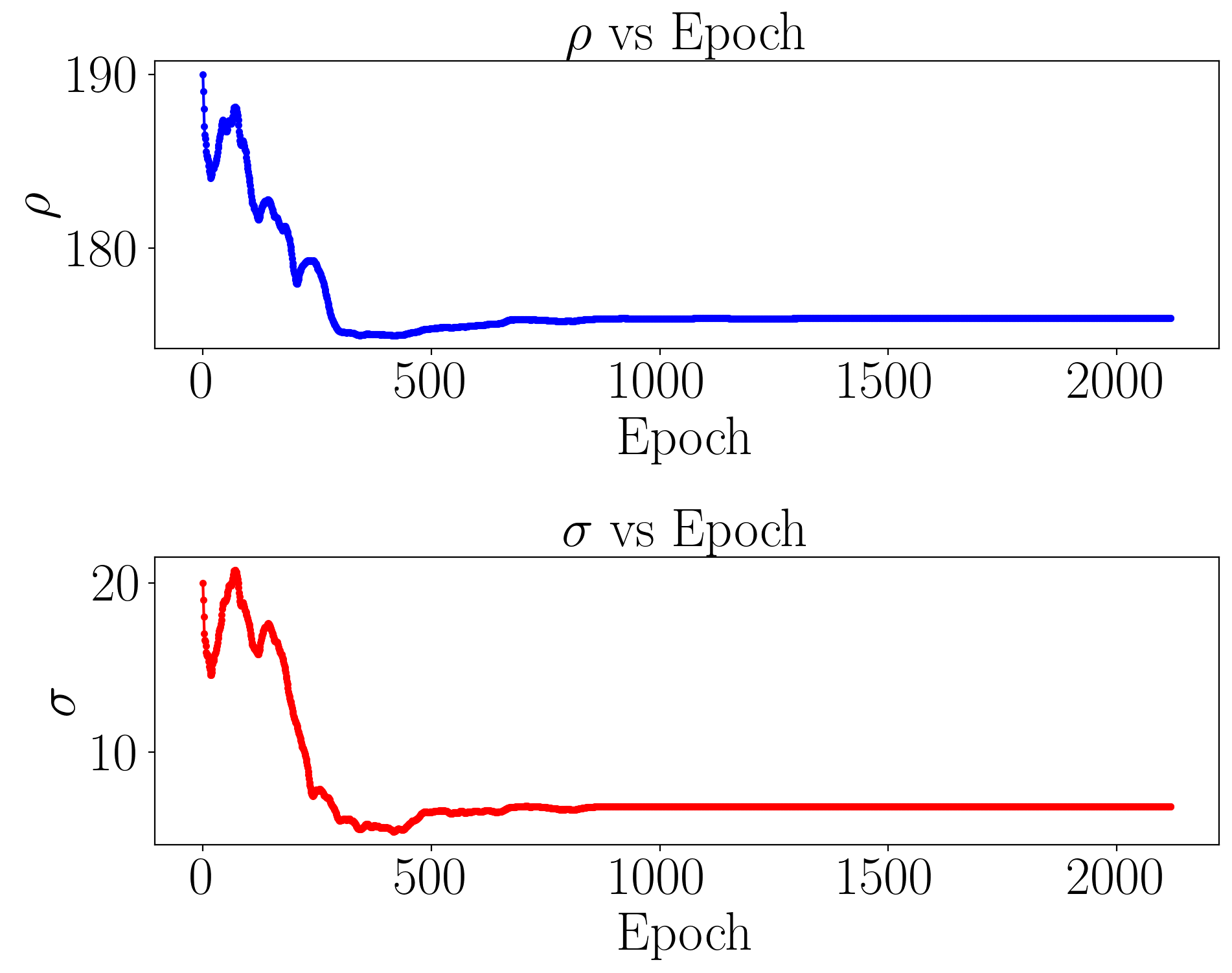}
  \caption{Lorenz system path components with respect to optimization epoch to demonstrate convergence.}
  \label{fig:lorenz_path_components}
\end{figure}

For the first test, we set the learning rate to 0.1 for the optimization. While this learning rate is too small to reach the periodic region in a reasonable amount of time, it is important to always start small with the learning rate to promote shorter paths. The resulting path without learning rate decay after 385 epochs is shown in Fig.~\ref{fig:lorenz0.1_final}.
We see that the path is very short in this case due to the small learning rate and it was never able to exit the chaotic region. In order to explore a larger region of the parameter space, the learning rate was then increased to one which is significantly larger than typical optimization problems, but the learning rate directly corresponds to the step size in the parameter space so it is justified in this problem as long as learning rate decay is used to ensure convergence. With a decay rate of 1\% per epoch, the resulting path after 382 epochs is shown in Fig.~\ref{fig:lorenz1.0_0.99_final}.
We see that the path was able to escape the chaotic region and reach a periodic solution, but the optimizer did not have enough momentum to reach the periodic solutions in the region for low values of $\sigma$. To allow for more exploration, in the final example the learning rate remained one and the decay rate was reduced to half a percent per epoch. In this case, the optimization was carried out to 2100 epochs to verify convergence, but the periodic solution was found after about 300 epochs. The resulting path and final trajectory are shown in Fig.~\ref{fig:lorenz1.0_0.995_final}
For completeness, the $\rho$ and $\sigma$ components of the path are also plotted with respect to epoch in Fig.~\ref{fig:lorenz_path_components} where it is clear that the path converges to a single pair of parameters at $\rho=175.996$ and $\sigma=6.776$.

\section{Conclusions}

Exploring dynamical system parameter spaces is a highly nontrivial task and there are many different approaches to solving the problem. We harnessed the connection between topology and dynamical systems to navigate parameter spaces by optimizing topological features of a dynamical system trajectory using persistence optimization. This resulted in the creation of a new language for defining topologically driven loss functions that map to different response characteristics. The loss function dictionary allows for promoting or avoiding limit cycles, fixed points and chaos for a dynamical system. Many choices need to be made by the user for this method to work such as the simulation time and sample frequency for simulations. If the transient behavior is not properly removed from the point cloud the persistence diagrams will be incorrect, but if the simulation time and sample frequency are too high, the computation times quickly become unreasonable. The user also needs to choose a learning rate and decay rate to optimize between exploration and exploitation. Once these factors have been tuned, the loss function can be intuitively designed to promote or avoid different features or regions of the parameter space. We believe that the most important decisions for the success of this method were adding gradient clipping to avoid exploding gradients in chaotic regions of the parameter space and balancing the loss function terms in the multi-objective optimization examples. These choices led to successful demonstrations over many examples from three different dynamical systems. In future work, it will be crucial to add functionality to account for unstable system responses, and avoid them at all costs or search for them for characterizing unsafe bounds on system parameters. In this work the parameter spaces were specifically chosen to not contain any unstable solutions. 

\section*{Code Availability}
We have made the code to reproduce this work publicly available in teaspoon \cite{Munch2018}, an open-source Python library for
Topological Signal Processing (TSP). Documentation and example code are also provided to aid with reproducibility of the results.

\section*{Acknowledgments}
This work is supported in part by Michigan State University and the National Science Foundation Research Traineeship Program (DGE-2152014) to Max Chumley.

\bibliographystyle{ieeetr}  
\bibliography{bibliography}

\begin{thebibliography}{10}

\bibitem{Kantz2004}
H.~Kantz and T.~Schreiber, {\em Nonlinear Time Series Analysis}.
\newblock Cambridge: Cambridge University Press, nov 2004.

\bibitem{Hirsch2003}
M.~W. Hirsch, S.~Smale, and R.~Devaney, {\em Differential Equations, Dynamical Systems, and an Introduction to Chaos (Pure and Applied Mathematics (Academic Press), 60.)}.
\newblock Academic Press, 2003.

\bibitem{Baker1996}
G.~L. Baker and J.~P. Gollub, {\em Chaotic Dynamics}.
\newblock Cambridge University Press, 1 1996.

\bibitem{Sayama2015}
H.~Sayama, {\em Introduction to the modeling and analysis of complex systems}.
\newblock Open SUNY Textbooks, 2015.

\bibitem{Kuznetsov1998}
Y.~A. Kuznetsov, {\em Elements of applied bifurcation theory}.
\newblock New York: Springer, 1998.

\bibitem{Seydel2009}
R.~U. Seydel, {\em Practical Bifurcation and Stability Analysis}.
\newblock Springer-Verlag GmbH, Nov. 2009.

\bibitem{Dankowicz2013}
H.~Dankowicz and F.~Schilder, {\em Recipes for Continuation}.
\newblock Society for Industrial and Applied Mathematics, 5 2013.

\bibitem{Sieber2007}
J.~Sieber and B.~Krauskopf, ``Control based bifurcation analysis for experiments,'' {\em Nonlinear Dynamics}, vol.~51, pp.~365--377, 2 2007.

\bibitem{Sieber2010}
J.~Sieber, B.~Krauskopf, D.~Wagg, S.~Neild, and A.~Gonzalez-Buelga, ``Control-based continuation of unstable periodic orbits,'' {\em Journal of Computational and Nonlinear Dynamics}, vol.~6, 9 2010.

\bibitem{Barton2011}
D.~A. Barton, B.~P. Mann, and S.~G. Burrow, ``Control-based continuation for investigating nonlinear experiments,'' {\em Journal of Vibration and Control}, vol.~18, pp.~509--520, 2 2011.

\bibitem{Bureau2013}
E.~Bureau, F.~Schilder, I.~F. Santos, J.~J. Thomsen, and J.~Starke, ``Experimental bifurcation analysis of an impact oscillator{\textemdash}tuning a non-invasive control scheme,'' {\em Journal of Sound and Vibration}, vol.~332, pp.~5883--5897, 10 2013.

\bibitem{Barton2013}
D.~A.~W. Barton and J.~Sieber, ``Systematic experimental exploration of bifurcations with noninvasive control,'' {\em Physical Review E}, vol.~87, p.~052916, 5 2013.

\bibitem{Barton2017}
D.~A. Barton, ``Control-based continuation: Bifurcation and stability analysis for physical experiments,'' {\em Mechanical Systems and Signal Processing}, vol.~84, pp.~54--64, 2 2017.

\bibitem{Godwin2017}
S.~Godwin, D.~Ward, E.~Pedone, M.~Homer, A.~G. Fletcher, and L.~Marucci, ``An extended model for culture-dependent heterogenous gene expression and proliferation dynamics in mouse embryonic stem cells,'' {\em npj Systems Biology and Applications}, vol.~3, 8 2017.

\bibitem{KRAUSKOPF2006}
B.~Krauskopf, H.~M. Osinga, E.~J. Doedel, M.~E. Henderson, J.~Guckenheimer, A.~Vladimirsky, M.~Dellnitz, and O.~Junge, ``A survey of methods for computing (un)stable manifolds of vector fields,'' in {\em World Scientific Series on Nonlinear Science Series B}, pp.~67--95, World Scientific, 3 2006.

\bibitem{Peeters2009}
M.~Peeters, R.~Vigui{\'{e}}, G.~S{\'{e}}randour, G.~Kerschen, and J.-C. Golinval, ``Nonlinear normal modes, part {II}: Toward a practical computation using numerical continuation techniques,'' {\em Mechanical Systems and Signal Processing}, vol.~23, pp.~195--216, 1 2009.

\bibitem{Huntley2017}
S.~Huntley, D.~Jones, and A.~Gaitonde, ``Bifurcation tracking for high reynolds number flow around an airfoil,'' {\em International Journal of Bifurcation and Chaos}, vol.~27, p.~1750061, 4 2017.

\bibitem{kuehn2008}
C.~Kuehn, ``Exploring parameter spaces in dynamical systems,'' 2008.

\bibitem{Hatcher2002}
A.~Hatcher, {\em Algebraic Topology}.
\newblock Cambridge University Press, 2002.

\bibitem{Kaczynski2004}
T.~Kaczynski, K.~Mischaikow, and M.~Mrozek, {\em Computational Homology}.
\newblock Springer, Jan. 2004.

\bibitem{Ghrist2008}
R.~Ghrist, ``Barcodes: The persistent topology of data,'' {\em Builletin of the American Mathematical Society}, vol.~45, pp.~61--75, 2008.
\newblock Survey.

\bibitem{Carlsson2009}
G.~Carlsson, ``Topology and data,'' {\em Bulletin of the American Mathematical Society}, vol.~46, pp.~255--308, 1 2009.
\newblock Survey.

\bibitem{Edelsbrunner2010}
H.~Edelsbrunner and J.~Harer, {\em Computational Topology: An Introduction}.
\newblock Rhode Island: American Mathematical Society, 2010.

\bibitem{Mischaikow2013}
K.~Mischaikow and V.~Nanda, ``Morse theory for filtrations and efficient computation of persistent homology,'' {\em Discrete \& Computational Geometry}, vol.~50, no.~2, pp.~330--353, 2013.

\bibitem{oudot2017persistence}
S.~Y. Oudot, {\em Persistence theory: from quiver representations to data analysis}, vol.~209 of {\em AMS Mathematical Surveys and Monographs}.
\newblock Rhode Island: American Mathematical Soc., 2017.

\bibitem{Munch2017}
E.~Munch, ``A user's guide to topological data analysis,'' {\em Journal of Learning Analytics}, vol.~4, pp.~47--61, jul 2017.

\bibitem{CohenSteiner2006}
D.~Cohen-Steiner, H.~Edelsbrunner, and J.~Harer, ``Stability of persistence diagrams,'' {\em Discrete {\&} Computational Geometry}, vol.~37, pp.~103--120, dec 2006.

\bibitem{Carriere2020}
M.~Carriere, F.~Chazal, M.~Glisse, Y.~Ike, and H.~Kannan, ``Optimizing persistent homology based functions,'' in {\em International conference on machine learning}, 10 2020.

\bibitem{Leygonie2021}
J.~Leygonie, S.~Oudot, and U.~Tillmann, ``A framework for differential calculus on persistence barcodes,'' {\em Foundations of Computational Mathematics}, vol.~22, pp.~1069--1131, 7 2021.

\bibitem{chen2019}
R.~T.~Q. Chen, Y.~Rubanova, J.~Bettencourt, and D.~Duvenaud, ``Neural ordinary differential equations,'' 2018.

\bibitem{Gameiro2016}
M.~Gameiro, Y.~Hiraoka, and I.~Obayashi, ``Continuation of point clouds via persistence diagrams,'' {\em Physica D: Nonlinear Phenomena}, vol.~334, pp.~118--132, 11 2016.

\bibitem{Atienza2020}
N.~Atienza, R.~Gonzalez-Díaz, and M.~Soriano-Trigueros, ``On the stability of persistent entropy and new summary functions for topological data analysis,'' {\em Pattern Recognition}, vol.~107, p.~107509, 2020.

\bibitem{Myers2019}
A.~Myers, E.~Munch, and F.~A. Khasawneh, ``Persistent homology of complex networks for dynamic state detection,'' {\em Physical Review E}, vol.~100, p.~022314, 8 2019.

\bibitem{Endres2018}
S.~C. Endres, C.~Sandrock, and W.~W. Focke, ``A simplicial homology algorithm for lipschitz optimisation,'' {\em Journal of Global Optimization}, vol.~72, pp.~181--217, 2018.

\bibitem{Genesio_2008}
R.~Genesio, G.~Innocenti, and F.~Gualdani, ``A global qualitative view of bifurcations and dynamics in the rössler system,'' {\em Physics Letters A}, vol.~372, pp.~1799--1809, Mar. 2008.

\bibitem{zhang2019}
J.~Zhang, T.~He, S.~Sra, and A.~Jadbabaie, ``Why gradient clipping accelerates training: A theoretical justification for adaptivity,'' 2019.

\bibitem{bischof2021}
R.~Bischof and M.~Kraus, ``Multi-objective loss balancing for physics-informed deep learning,'' 2021.

\bibitem{Xiao2024}
B.~Xiao, ``Strategies for balancing multiple loss functions in deep learning.'' Medium, 2024.

\bibitem{Myers_2020_pendulum}
A.~Myers and F.~A. Khasawneh, ``Dynamic state analysis of a driven magnetic pendulum using ordinal partition networks and topological data analysis,'' in {\em Volume 7: 32nd Conference on Mechanical Vibration and Noise ({VIB})}, American Society of Mechanical Engineers, aug 2020.

\bibitem{Munch2018}
E.~Munch, ``Teaspoon.'' https://github.com/lizliz/teaspoon, 2018.

\end{thebibliography}

\end{document}